\documentclass[11pt]{ip-journal}

\usepackage[new]{old-arrows}
\usepackage{longtable}
\usepackage{booktabs}
\usepackage{array}
\usepackage{multirow}
\usepackage{supertabular}
\usepackage{amssymb}
\usepackage{amsmath}
\usepackage{mathrsfs}
\usepackage{titletoc}
\usepackage{graphicx}
\usepackage{latexsym}
\usepackage[alphabetic]{amsrefs}

 \newtheorem{thm}{Theorem}[section]

 \newtheorem{defn}[thm]{Definition}
 \newtheorem{rem}[thm]{Remark}
  \newtheorem{thm*}{Theorem}
\newtheorem{lem*}{Lemma}
\newtheorem{cor*}{Corollary}
\newtheorem{prop*}{Proposition}

\begin{document}

\title[Minimizing cones associated with isoparametric foliations\ \ ]{Minimizing cones\\ associated with isoparametric foliations\ \ }
\author{Zizhou Tang and Yongsheng Zhang}
\address{
Chern Institute of Mathematics \& LPMC, Nankai University, Tianjin 300071, P. R. China
}\email{zztang@nankai.edu.cn}
\address{
{Current Address:\ School of Mathematical Sciences \& IAS, Tongji University, Shanghai, 200092,
 P. R. China}
}
\email{yongsheng.chang@gmail.com}
\date{\today}

\begin{abstract}

{
Associated with isoparametric foliations of unit spheres,
there are two classes of minimal surfaces $-$ minimal isoparametric hypersurfaces and focal submanifolds.
By virtue of their rich structures,
we find new series of minimizing cones.
They are cones over focal submanifolds 
and
cones over suitable products among
these two classes.
Except in low dimensions,
all such cones are shown minimizing.
}
     \end{abstract}
     

\maketitle
\titlecontents{section}[0em]{}{\hspace{.5em}}{}{\titlerule*[1pc]{.}\contentspage}
\titlecontents{subsection}[1.5em]{}{\hspace{.5em}}{}{\titlerule*[1pc]{.}\contentspage}
{\setcounter{tocdepth}{2} \small \tableofcontents}
\section{Introduction}\label{Section1}
          Minimizing hypercones have fascinating originations and active developments in geometric measure theory, for example,
          the pioneering works:
           Fleming \cite{Fleming},
          De Giorgi \cite{DG}, 
          Almgren \cite{Almgren},
          Simons \cite{Simons},
          Bombieri-De Giorgi-Giusti \cite{BDG} for solving the celebrated Bernstein problem.

         Many types of homogeneous minimizing hypercones beyond Simons cones were discovered
         by Lawson in \cite{BL}.
         Later \cite{PS1, PS2} added that $C_{2,4}$ over minimal Clifford torus
         $S^{2}\left(\sqrt{\frac{1}{3}}\right)
         \times 
         S^{4}\left(\sqrt{\frac{2}{3}}\right)
         \subset S^7(1)$ is minimizing  (also see \cite{Z2})
         whereas $C_{1,5}$ is only stably minimal.
         Although lots of examples had been found,
         a complete  list of homogeneous minimizing hypercones remained unknown until the important work 
         of Lawlor \cite{Law}.
         Rather than in the entire Euclidean space,
         Lawlor searched for a structure similar to the characterization foliation of a minimizing hypercone
         discovered 
         by Hardt-Simon \cite{HS}
         in certain angular neighborhood of the cone under consideration.

         More flexible than homogeneous foliations are the isoparametric ones.
         In 1985, Ferus and Karcher \cite{FK}
         successfully {constructed} the characterization foliation 
           for almost all minimal isoparametric hypercones of OT-FKM type (see \S \ref{Section2} below or \cite{OT, FKM} for more details), and hence showed each of them minimizing.%
Thereafter, to figure out a complete classification of minimizing isoparametric hypercones became important in the subject.
         Based on the method of \cite{FK}, an explicit list was given in \cite{Wang}.
However, the arguments for strict minimality in his \S 9 were invalid
         due to the inaccurate Remark 9.17 on a correct formula in \cite{HS}.
         We shall give an alternative proof in \S 2.
        
          Since a tangent cone (at some point) of an area-minimizing rectifiable current is itself area-minimizing
              (Theorem 5.4.3 in \cite{F}, also see Theorem 35.1 and Remark 34.6 (2) in \cite{LS}),
              the study of minimizing cones of higher codimensions is of equal importance.
              Besides minimal isoparametric hypersurfaces,
                another class of minimal submanifolds associated with isoparametric foliations are focal submanifolds.
                By the efforts of \cite{Law, Kerckhove, HKT, Kanno, OS},
                cones over focal submanifolds of homogeneous isoparametric foliations 
                with $g=3,4,6$ distinct principal curvatures and multiplicities $(m_1,m_2)\neq (1,1)$
                had been shown minimizing.
                %
              In this paper, we shall consider cones over general focal submanifolds.
             
              {\ }
              
              In \S\ref{Section3},
%
         we establish
        \\{\ }
        
        \textbf{Theorem 1.}\label{thm1}
       {\it  Cones over focal submanifolds of isoparametric foliations
       with $g=4$ 
       and 
       $(m_1,m_2)\neq (1,1)$ are area-minimizing.}

        \begin{rem}
       Lots of them are inhomogeneous v.s. those equivariant minimizing cones constructed in \cite{XYZ2}.
       {It should be remarked that for some inhomogeneous isoparametric foliations their focal submanifolds can possibly be homogeneous.
       For details readers are referred to \cite{FKM}.
       }
       \end{rem}

       \begin{rem}
       For $g=4, m_1=m_2=2$, both $C(M_ +)$ and $C(M_-)$, the cones over the focal submanifolds $M_+$ and $M_-$,  are minimizing and of dimension $7$ in $\mathbb R^{10}$.
       Moreover,
                     since the vanishing angle is less than $\frac{\pi}{8}$ $=\frac{1}{2}\cdot\frac{\pi}{4}$,
                     their union $C(M_+\bigsqcup M_-)$ wth either orientation combination is minimizing as well.
       \end{rem}

        \begin{rem}
       When $M_t$ for $t=+$ or $-$ is nonorientable, $C(M_t)$ is minimizing in the sense of mod $2$ (see \cite{Ziemer}).
         For $g=4$ and $(m_1,\, m_2)=(1,k)$, $M_+$ is diffeomorphic to $V_2(\mathbb R^{k+2})$ and $M_-$ is isometric to $S^1\times S^{k+1}/\mathbb Z_2$ (see \cite{TY}).
       It can be proved without difficulty that the latter is orientable if and only if $k$ is even.
       \end{rem}
         {Our investigations} for $g=3\text{ and }6$ 
                 {lead} to
       \\{\ }
       
       \textbf{Theorem 2.}\label{thm2} (\cite{Kerckhove, Kanno, OS})
       {\it  Cones over focal submanifolds of isoparametric foliations 
       with {$g=3,6$} and $(m_1,m_2)\neq (1,1)$ are area-minimizing.}
       \\{\ }
       
        \textbf{Question.}
        It will be interesting, for a complete classification, to ask
        whether cones over focal submanifolds of isoparametric foliations with $g=3,4,6$ and $(m_1,m_2)=(1,1)$
        are minimizing.
        By Cartan, for $g=3$ and $(m_1,m_2)=(1,1)$, $M_\pm$ are both isometric to $\mathbb RP^2$ with 
          constant Gaussian curvature;
          for $g=6$ and $(m_1,m_2)=(1,1)$, $M_\pm$ are diffeomorphic to $S^3\times \mathbb RP^2$ (see, for example, \cite{Miyaoka}).
           \\{\ }

                         {Further, we consider the minimal products of minimal submanifolds in \S\ref{Section4}
                         and gain
                         }
       \\{\ }
       
        \textbf{Theorem 3.}\label{thm3}
       {\it  Cones over the minimal products of 
       focal submanifolds of isoparametric foliations
       with {$g=3,4,6$} and $(m_1,m_2)\neq (1,1)$ 
       are area-minimizing.}
\\{\ }

       \textbf{Theorem 4.}\label{thm4}
       {\it  Cones, of dimension no less than $10$, over the minimal products
       of 
       focal submanifolds
       are area-minimizing.}
      \\
      

                We continue to 
                study cones over
       products among focal submanifolds and minimal isoparametric hypersurfaces.
                {In particular, we show
       \\{\ }

        \textbf{Theorem 5.}\label{thm5}
       {\it  Cones, of dimension no less than $37$, over the minimal products
       among focal submanifolds
       and minimal isoparametric hypersurfaces
       are area-minimizing.}}
\\{\ }

       The paper is organized as follows.
       In \S 2 we review
       some fundamental results related to the topic of our paper.
       In \S3.1 we briefly go through the ideas of Lawlor's curvature criterion
       and apply it to cones over focal submanifolds in \S3.2 and \S3.3 for $g=4$ and $g=3$, $6$ respectively.    
       In \S 4, we first consider the minimal products of two focal submanfolds with $g=3,4,6$, 
       then the multiple case based on the observation Remark \ref{rknice},
       and finally include 
       $g=2$ 
       with the aid of Lawlor's results on cones over products
       purely of spheres.
       Extending \S\ref{Section4}, the last section is devoted to the case of products
       among minimal isoparametric hypersurfaces and focal submanifolds.
       Our arguments heavily rely on the foliation structures.


 \section{Preliminaries}\label{Section2}

              A closed (embedded) hypersurface $M$ in 
the
              unit sphere $S^{n-1}
              \subset \mathbb R^n$ is called isoparametric, by E. Cartan, if it has constant principal curvatures.
             When $M$ is isoparametric, so are its parallel hypersurfaces.
             In this way, a foliation of hypersurfaces appears with two exceptional leaves $-$ focal submanifolds of higher codimensions.
             It is well known that
             there is one and only one minimal hypersurface among parallel isoparametric hypersurfaces
             of an isoparametric foliation.

          Let $\xi$ be a unit normal vector field along $M$,
               $g$ the number of distinct principal curvatures of $M$,
          $\cot \theta_\alpha$ $(\alpha=1,\cdots, g; 0<\theta_1<\cdots<\theta_g<\pi)$
           the principal curvatures with respect to $\xi$
           and
           $m_\alpha$ the multiplicity of $\cot\theta_\alpha$.
           By a purely topological method, M\"{u}nzner proved 
           an elegant result that $g$ must be $1,2,3,4$ or $6$
           and $m_\alpha = m_{\alpha+2}$ (indices mod $g$).
           According to values of $g$, we have the followings.
           \\{ \ }

           $g=1$. The foliation is trivial,
                 namely given by level sets of a height function restricted to $S^{n-1}
                $
                 with respect to a {nonzero} direction {in $\mathbb R^n$}.

           $g=2$. Regular leaves are of type $S^p\times S^{n-2-p}$ and
           the corresponding focal submanifolds are great spheres $S^p$ and $S^{n-2-p}$.

           $g=3$. $m_1= m_2 =m_3$ has to take values among $1, 2, 4$ and $8$.
           Cartan
            \cite{Ca}
           showed that these are homogeneous foliations. 
                      Moreover,
                      all the isoparametric hypersurfaces 
                      are precisely
                      the tubes of constant radius over the standard Veronese embedding of
                      $\mathbb FP^2$ for $\mathbb F=\mathbb R$, $\mathbb C$, $\mathbb H$ (quaternions),
                      $\mathbb O$ (Cayley numbers) in $S^4$, $S^7$, $S^{13}$, $S^{25}$, respectively.

           $g=6$. $m_1=m_2=m$
            by M\"unzner \cite{Munzner} and
           $m$ has to be $1$ or $2$ by Abresch \cite{Abresch}.
           For $m=1$, Dorfmeister and Neher \cite{DN} proved that the foliation is homogeneous.
           Very recently, Miyaoka \cite{Miyaoka, Miyaoka2}
           established the same conclusion for $m=2$.

           $g=4$.
           By recent beautiful results of Cecil-Chi-Jensen \cite{CCJ} and Chi \cite{Chi},
           the classification {finally} gets complete.
           Such a foliation must be either of OT-FKM type or homogeneous with $(m_1,m_2)=(2,2)$ or $(4,5)$.
%
           Here the OT-FKM type means that
           leaves of the foliation are given
           by
           \[
           M_t=\bold F^{-1}(t)\bigcap S^{n-1} \ \ \ \text{ for } t\in[-1,+1]
           \]
           where
            \[
            \bold F(x)=\; <x,x>^2-2\sum_{i=0}^{m}<P_ix,x>^2,\ \ \ m\geq 1
            \]
            for self-adjoint endomorphisims $P_i:\mathbb R^n\rightarrow \mathbb R^n$,
            for $i=0,\cdots,m$,
            with relation $P_iP_j+P_jP_i=2\delta_{ij}\text{Id}$.
            Such structure exists only when $n=2l$ and $l=k\delta(m)$ (see \cite{FKM}).
            Its multiplicities satisfy $(m_1,m_2)=(m, l-m-1)$.
            \\{\ }

           The first inhomogeneous foliation of OT-FKM type occurs when $(m_1,m_2)=(3,4)$.
            For all other inhomogeneous foliations, $3\leq m_1<m_2$ and $m_1+m_2\geq 11$.
            Theorem 2 of \cite{FK} states that
                         each inhomogeneous minimal isoparametric hypersurface $M^*$ of OT-FKM type
                         {spans} a minimizing hypercone
                   \[
                   C(M^*)=\{tx: x\in M^*\text{ and } t\in [0,\infty)\}.
                   \]
             In fact their argument works for all minimal isoparametric hypersurface of OT-FKM type with $m_1+m_2\geq 11$.
         \\{\ }

                { Combined with} the {above} classification theorem of isoparametric foliations,
            Theorem 2 of \cite{FK}
             and 
             classification on homogeneous minimizing hypercones 
            confirm the classification in \cite{Wang}.
            
              {\ }
              
            \textbf{Theorem.} ({\cite{Wang})}
            Let $M^*$ be a minimal isoparametric hypersurface in $S^{n-1}$
            with $g\geq 2$.
            Then $C(M^*)$ is minimizing if and only if $n\geq 4g$ and $(g,m_1,m_2)\neq (2,1,5)$ or $(4,1,6)$.

            {\ }

            Here is a table of 
            homogeneous minimal hypercones about their being minimizing, stable minimal, and unstable minimal,            
            according to values of $(g, m_1, m_2)$
            (cf. \cites{Law, MO, Z2}).
            
  \makeatletter
\def\hlinew#1{%
  \noalign{\ifnum0=`}\fi\hrule \@height #1 \futurelet
   \reserved@a\@xhline}
\begin{center}
\begin{longtable}[h]{cp{2.2cm}<{\centering}p{2cm}<{\centering}p{1.9cm}<{\centering}p{1.9cm}<{\centering}}
\specialrule{0em}{3pt}{3pt}          
\hlinew{0.8pt}
\multirow{2}{*}{{\small $(g,m_1,m_2)$}}       &          \multicolumn{1}{c|}{ \multirow{2}{*}{\small Type of $M^*$}}
                                                                                &                                                                                 
                                                                           \multicolumn{3}{c}{\small Strongest minimality of $C(M^*)$ among}\\
                                                                                \cline{3-5}
                                                                                \multicolumn{2}{c|}{ }
                                                                            &  \small {Minimizing} &  {\small Stable} &  { {\small Unstable}}\\\hlinew{0.8pt}
\endfirsthead
\multicolumn{5}{l}%
{
\textit{\tiny Continued from previous page}}\\
\hlinew{0.8pt}
\multirow{2}{*}{{\small $(g,m_1,m_2)$}}       &          \multicolumn{1}{c|}{ \multirow{2}{*}{\centering{\small Type of   $M^*$}}}
                                                                                &   
                                                                                
                                                                                \multicolumn{3}{c}{ {\small Strongest minimality of $C(M^*)$} among}\\
                                                                                \cline{3-5}
                                                                                \multicolumn{2}{c|}{}
                                                                            & {\small {Minimizing}} &  {\small Stable} &   { {\small Unstable}}\\\hlinew{0.8pt}
\endhead
\specialrule{0em}{3.5pt}{3.5pt} 
\multicolumn{5}{r}{\textit{\tiny Continued on next page}}\\
\endfoot
\endlastfoot

\specialrule{0em}{3.5pt}{.5pt} 
 {\footnotesize $g=1$}                    &                                                         {\footnotesize $S^{n-2}$}        & \checkmark & &\\
\specialrule{0em}{3.5pt}{3.5pt}          
{\footnotesize $(2,p,n-2-p)$}     &    {\scriptsize $S^p\times S^{n-2-p}$}     &           {\tiny  $n=8,p\neq1,5$; $n>8$}
                                                 &                     {\tiny $n=8,p=1,5$}          &                 {\scriptsize   $n<8$}   \\
\specialrule{0em}{3.5pt}{3.5pt}                                                           
{\footnotesize $(4,1,p-2)$}        &     {\small   $\frac{SO(p)\times SO(2)}{SO(p-2)\times \mathbb Z_2}$}
   &             {\scriptsize {$p>8$}}   &                         {\scriptsize   $p=8$}              &          {\scriptsize  $3\leq p\leq7$}           \\
\specialrule{0em}{3.5pt}{3.5pt}              
{\footnotesize $(4,2,2p-3)$}        &     {\small   $\frac{S\left(U(p)\times U(2)\right)}{SU(p-2)\times T^2}$}      &            {\scriptsize  $p\geq4$}                   &                                 &        {\scriptsize $p=2,3$}          \\
\specialrule{0em}{3.5pt}{3.5pt}           
{\footnotesize $(4,4,4p-5)$}        &     {\small   $\frac{Sp(p)\times Sp(2)}{Sp(p-2)\times Sp(1)^2}$}    &            {\scriptsize   $ p\geq2$}               &                                 &                 \\
\specialrule{0em}{3.5pt}{3.5pt}          
{\footnotesize $(4,4,5)$}       &                        {\small   $\frac{U(5)}{SU(2)\times SU(2)\times T^1}$}    &         \checkmark                   &                &                   \\
\specialrule{0em}{3.5pt}{3.5pt}           
{\footnotesize $(4,6,9)$}   &   {\small   $\frac{Spin(10)\cdot T}{SU(4)\cdot T}$}      &         \checkmark                &                     &                   \\
 \specialrule{0em}{3.5pt}{3.5pt}          
{\footnotesize $(3,1,1)$}     &                    {\small   $\frac{SO(3)}{\mathbb Z_2+ \mathbb Z_2}$}   &                                      &                                &     \checkmark           \\
 \specialrule{0em}{3.5pt}{3.5pt}          
{\footnotesize $(3,2,2)$}     &                        {\small    $\frac{SU(3)}{T^2}$}                                    &                              &                                       &     \checkmark            \\
 \specialrule{0em}{3.5pt}{3.5pt}          
 {\footnotesize $(3,4,4)$}    &              {\small                 $\frac{Sp(3)}{Sp(1)^3}$}                              &      \checkmark                     &                                        &                \\
 \specialrule{0em}{3.5pt}{3.5pt}          
{\footnotesize $(3,8,8)$}   &                   {\small            $\frac{F_4}{Spin(8)}$}                                      &       \checkmark                  &                                      &                 \\
\specialrule{0em}{3.5pt}{3.5pt}           
 {\footnotesize $(4,2,2)$}     &                      {\small             $\frac{SO(5)}{T^2}$}                                       &                               &                                      &    \checkmark          \\
 \specialrule{0em}{3.5pt}{3.5pt}          
 {\footnotesize $(6,2,2)$}   &                        {\small         $\frac{G_2}{T^2}$}                                    &                                      &                            &   \checkmark            \\
 \specialrule{0em}{3.5pt}{3.5pt}          
{\footnotesize $(6,1,1)$}   &  {\small   $\frac{SO(4)}{\mathbb Z_2+ \mathbb Z_2}$}                 &     &                                &  \checkmark              \\
\specialrule{0em}{3pt}{2pt}
\hlinew{0.8pt} 
\specialrule{0em}{3pt}{1pt} 
\end{longtable}
\end{center}

          { For a minimal isoparametric hypersurface $M^*$,
            b}y
           \S 4.3 of \cite{Law},
            Lawlor's criterion about $C(M^*)$
                   is not only sufficient but also necessary
              {for its being minimizing}.
           What is more, whenever the criterion {applies},
           the cone becomes strictly area-minimizing automatically.
           Hence the next statement follows.
\\{\ }

            \textbf{Proposition.} {(\cite{Wang})}
            Let $M^*$ be a minimal isoparametric hypersurface in $S^{n-1}$.
            {Then} $C(M^*)$ is minimizing
             {if and only if}
            it is strictly area-minimizing.


\section{Cones over focal submanifolds}\label{Section3}

       It is well known that
       focal submanifolds of isoparametric foliations are minimal submanifolds {in spheres}.
       They generate minimal cones.
       Based on their special second fundamental forms,
       we shall {employ} the curvature criterion of Lawlor.

       To be self-contained, we briefly review Lawlor's method in \S\ref{LCC}.
       A proof of Theorem \ref{thm1} will be given in \S\ref{g4}.
       Along the same line, homogeneous cases for $g=3$ and $6$ will be discussed in \S\ref{g36}.

\subsection{Curvature criterion of Lawlor}\label{LCC}

       Let us {mention} some notations and concepts.
       \begin{defn}
       Let $B$ be a submanifold of dimension $k-1$ in $S^{n-1}$ and $C=C(B)$.
       Fix $p\in B$.
       Let $N^{n-k}$ be the $(n-k)$-dimensional great sphere which intersects $B$ at $p$ orthogonally.
       For $0<\eta<\pi$, denote by $U_p(\eta)$ the open $\eta$-disk centered at $p$ in $N^{n-k}$.
       Then $W_p(\eta)=C(U_p(\eta))$ is called
                    the $\mathrm{\mathbf{\eta\text{-}normal\ wedge}}$
                     through $p$.
       Conventionally, we leave out the origin, so that one can talk about ``nonintersecting normal wedges".
       We name $\bigcup_{p\in B}W_p(\eta)$ the
       $\mathrm{\mathbf{\eta\text{-}angle\ neighborhood}}$ of $C$.
       The $\mathrm{\mathbf{normal\ radius}}$
       of $C$ at $p$ means the largest angular radius $\eta$ for $W_p(\eta)$
        intersecting $C$ only in the ray $\overrightarrow{0p}$.
       \end{defn}
       \begin{figure}[h]
                        \begin{center}
                        \includegraphics[scale=0.5]{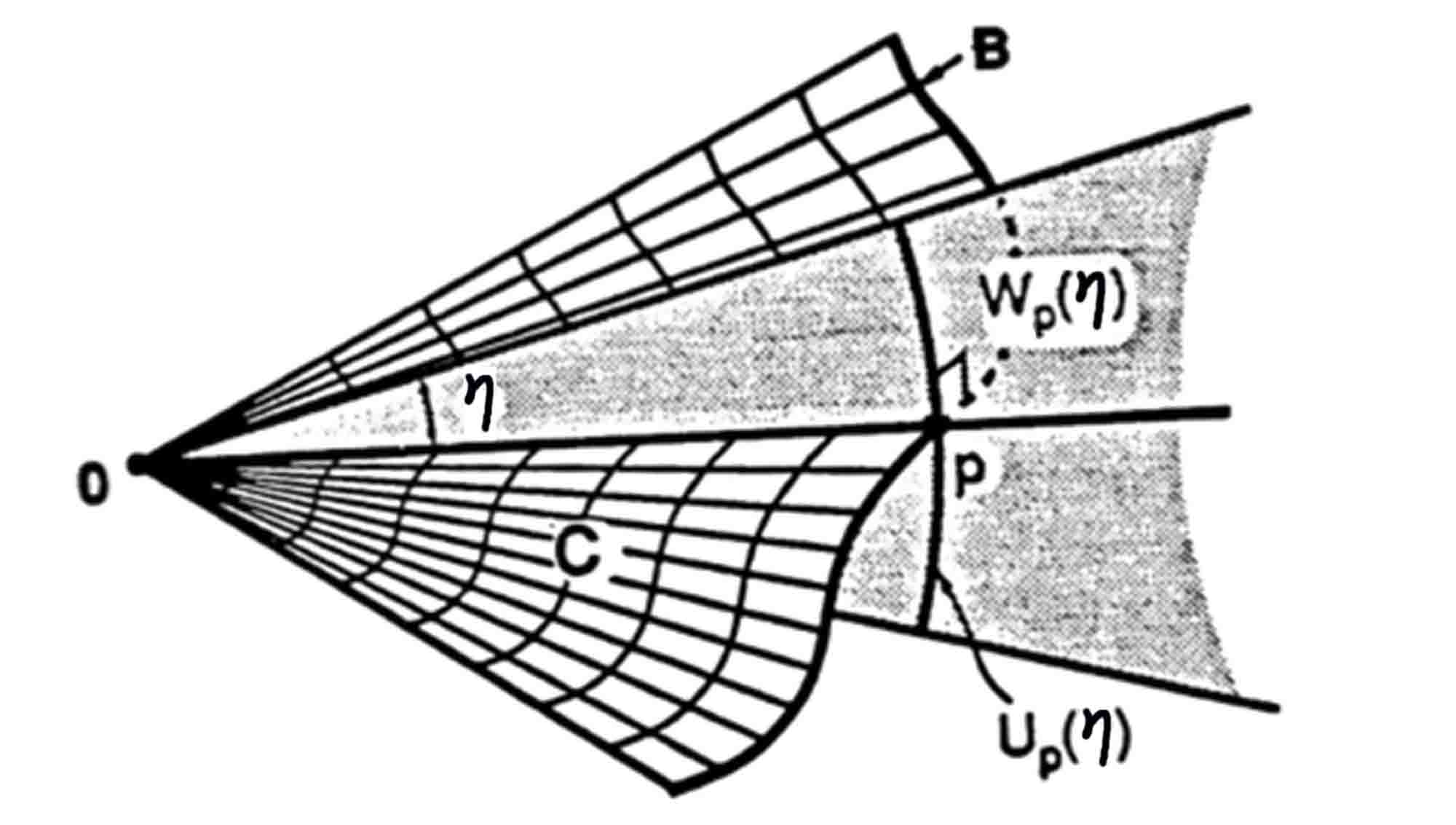}
                        {\\ \small $\eta$-normal wedges of a cone}
                        \end{center}
         \end{figure}

          Since the characterization foliation for an area-minimizing hypercone in \cite{HS} has the property of homothety,
             so do its perpendicular integral curves (outside the origin).
             The latter naturally induces an area-nonincreasing projection {for hypersurfaces} to the cone.
             Instead, Lawlor looked for similar structure in some $\eta$-angular neighborhood of a minimal cone $C$
                 rather than in $\mathbb R^n$.
             If the boundary of the neighborhood happens to be mapped to the origin under the projection,
             then one can send everything outside the neighborhood to the origin.
             In this way an area-nonincreasing projection can be produced.
{For higher codimensions,}
             Lawlor considered the structure given by rotation of a suitable curve $\gamma_p$
             in each normal wedge $W_p(\eta)$ for $p\in B$.
             Positive homotheties of the surface rotated by $\gamma_p$ are required to foliate the normal wedge.
             Under the assumption that
                 elements of $\{W_p(\eta) : p\in B\}$ do not intersect each other,
             the existence of preferred area-nonincreasing projection in the $\eta$-angular neighborhood
             is equivalent to
             that
               the {ordinary differential inequality}:
                          \begin{equation*}\label{cineq0}
                           \left(h(t)-\frac{t}{k}h'(t)\right)^2+\left(\frac{h'(t)}{k}\right)^2\leq \Big(\det(I-th_{ij}^{\nu})\Big)^2,\ \ \ \ \ \ h(0)=1,
                          \end{equation*}
               where $t=\tan(\theta)$ for $\theta\in [0,\eta)$
               and
{the $(k-1)\times(k-1)$ matrix
                $(h_{ij}^\nu)$} means the second fundamental form 
{of $B$
               at $p$ for a unit normal $\nu$ in $S^{n-1}$},
                has a solution
                which can reach zero
                for each $\nu$.
                Letting $\nu$ range over all unit normals at $p$,
                the requirement becomes that
                 \begin{equation}\label{cineq}
                           \left(h(t)-\frac{t}{k}h'(t)\right)^2+\left(\frac{h'(t)}{k}\right)^2\leq \Big(q(t)\Big)^2,\ \ \ \ \ \ h(0)=1,
                          \end{equation}
or equivalently,
               \begin{equation}\label{ncineq}
               \begin{split}
             & \dfrac{k}{t^2+1}\left(th-\sqrt{(t^2+1)q^2(t)-h^2}\right)\leq h'(t)\\
             &{\ \ \ \ \ \ \ \ \ \ \ \ \ \ \ \ \ \ \ }
              \leq \dfrac{k}{t^2+1}\left(th+\sqrt{(t^2+1)q^2(t)-h^2}\right),\ \ \ \ h(0)=1,
              \end{split}
                          \end{equation}     
             where $q(t)=\inf_\nu \det(I-th_{ij}^{\nu})=1+q_2t^2+\cdots$, 
             supports a solution which reaches zero.
         {It can be seen} that the solution $h_0(t)$
                  which attains zero fastest, if exists, must satisfy 
\begin{equation}\label{h0}
                 h_0'(t)= \dfrac{k}{t^2+1}\left(th_0-\sqrt{(t^2+1)q^2(t)-h_0^2}\right),\ \ \ \ \ \ h_0(0)=1.
                  \end{equation}
             \begin{defn}
             Suppose $h_0(t)$ gets zero at $t=\tan(\theta_{0}(p))$.
             We call
$\theta_0(p)$ the vanishing angle at $p$ and
         $\theta_0=\max_{p\in B}\theta_{0}(p)$ the vanishing angle of $C$.
             \end{defn}
             $\theta_0$ stands for the narrowest (uniform) size for
               the preferred area-nonincreasing projection.      
         {Hence, there will be} two things to check     
         {to apply} Lawlor's criterion:
            {\ }

             1. There exists a finite vanishing angle $\theta_0$.
             Here we would like to remark that 
{in general}
             \eqref{cineq} and \eqref{h0} may support no solutions which can touch zero.
            {\ }

             2. The $\theta_0$-normal wedges do not intersect.
             \\
             {\ }

             A technical point is how to control $q(t)$ in practice.
             Corollary 1.3.3 in \cite{Law} says that,
                        for $\mathscr M=
                        {k-1}\geq 2$
                         and $t\in[0,\frac{1}{\alpha}\sqrt{\frac{\mathscr M}{\mathscr M-1}}]$,
                        \begin{equation}\label{str}
                        \det (I-th^\nu_{ij})\geq F(\alpha, t, \mathscr M)
                        \end{equation}
               where
\begin{eqnarray}
                          {\ \ \ \ \ \ } F(\alpha, t, \mathscr M)&=&\left(1-\alpha t\sqrt{\frac{\mathscr M-1}{\mathscr M}}\right)\left(1+\frac{\alpha t}{\sqrt{\mathscr M(\mathscr M-1)}}\right)^{\mathscr M-1} ,
                          \label{Fofstr}
                          \\
                             \specialrule{0em}{2pt}{2pt} 
               \alpha&=&\sup_\nu\|(h_{ij}^\nu)\| \;=\; \sup_\nu \sqrt{\ \sum_{i,j}(h^\nu_{ij})^2}. \label{aofstr}
\end{eqnarray}
               The inequality is sharp; for certain matrices, equality holds for all $t$.
               (For example, for the classical coassociative Lawson-Osserman cone in $\mathbb R^7$, see \cite{XYZ2}.)
               Note that
               {the expression}
               $F(\alpha, t, \mathscr M)$ is nonincreasing in $\mathscr M$.
               Limiting $\mathscr M\rightarrow \infty$ leads to Corollary 1.3.4
                       \begin{equation}\label{nstr}
                       \det (I-th^\nu_{ij})\geq  F(\alpha, t,
{k-1}
                       )> (1-\alpha t)e^{\alpha t}.
                        \end{equation}

                 Based on \eqref{str} and \eqref{nstr}, Lawlor
                 considered
                 \begin{eqnarray}
                       &\ \ \ h'(t)= \dfrac{k}{t^2+1}\left(th-\sqrt{(t^2+1)(F(\alpha, t,
{k-1}
                        ))^2-h^2}\right),   \ &h(0)=1;\label{streq}\\
                           \specialrule{0em}{2pt}{2pt} 
                        &h'(t)= \dfrac{k}{t^2+1}\left(th-\sqrt{(t^2+1)((1-\alpha t)e^{\alpha t})^2-h^2}\right), \ &h(0)=1,\label{nstreq}
                       \end{eqnarray}
                 respectively.
                            Let $\theta_F(p)$ and $\theta_c(p)$ be the corresponding vanishing angles.
                 Then
                     \[
                     \theta_0(p)\leq \theta_F(p)<\theta_c(p).
                     \]
                 Lawlor
                 gained
                 the following table of upper bounds of vanishing angles for $\dim (C)$ and $\alpha^2$.
                  \begin{figure}[ht]
                              \begin{minipage}[c]{0.53\textwidth}
                              \includegraphics[scale=0.285]{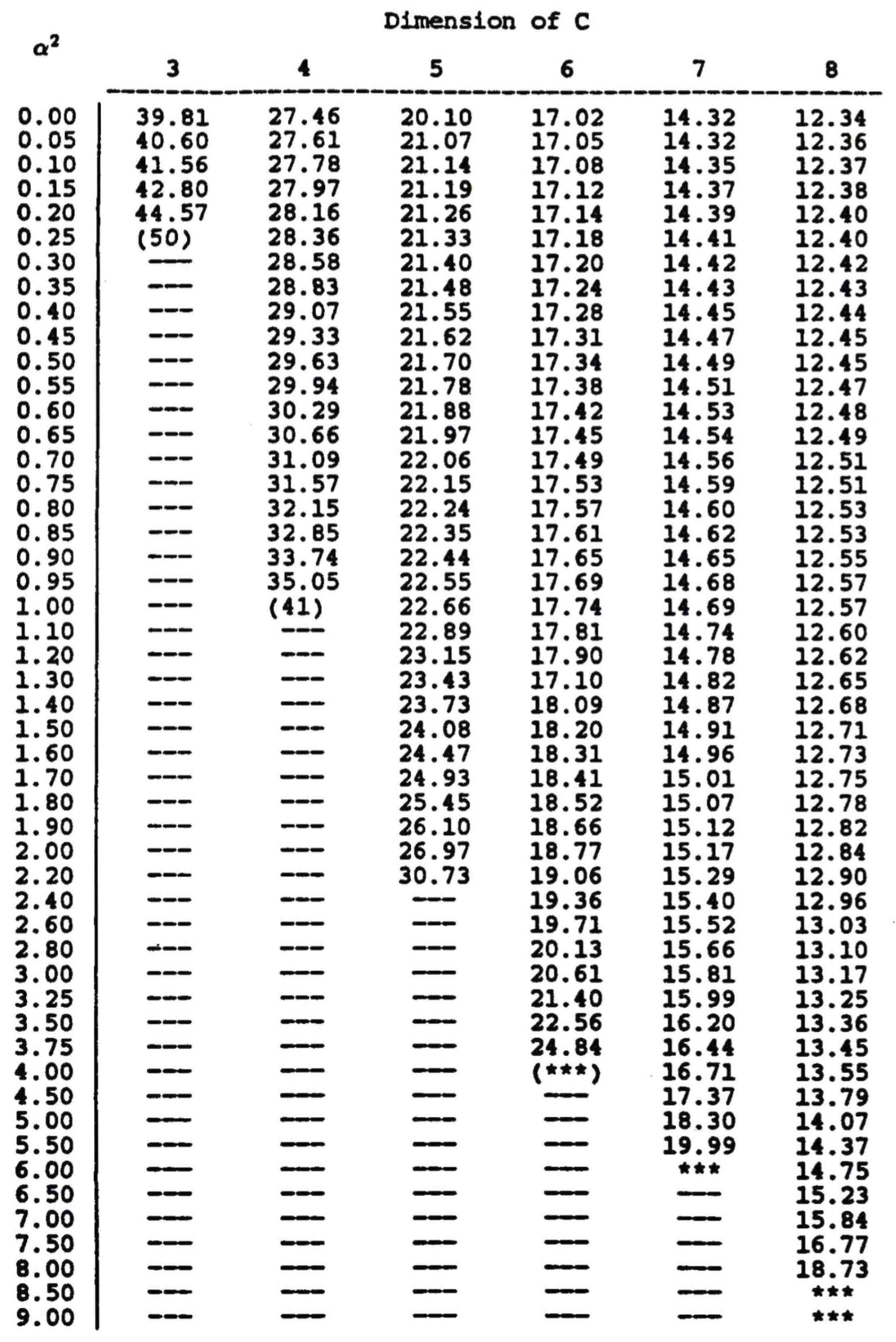}
                              \end{minipage}%
                          \begin{minipage}[c]{0.5\textwidth}
                           \includegraphics[scale=0.315]{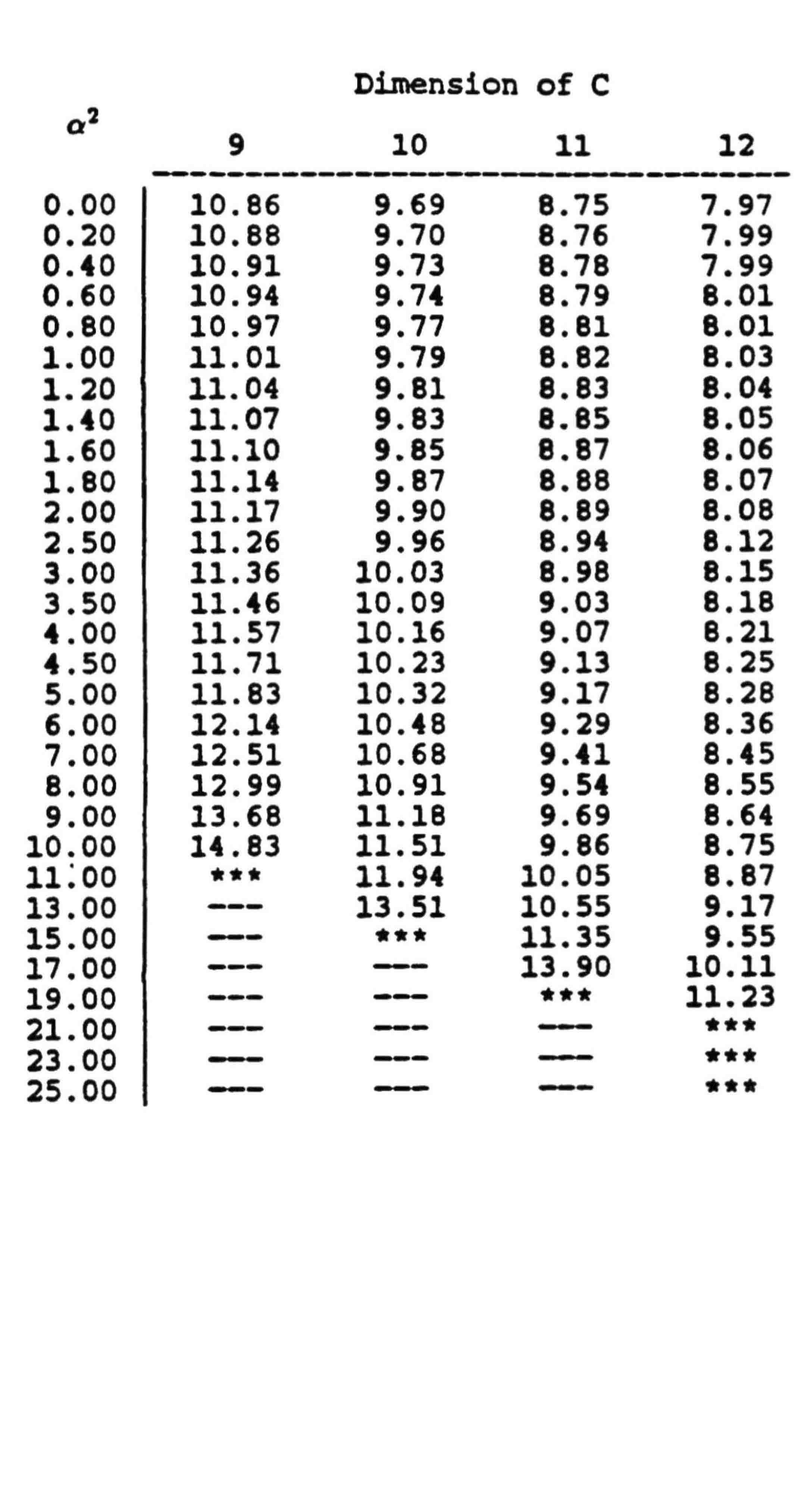}
                           \end{minipage}
                           {\\ \small {Table of upper bounds of vanishing angles}}
                    \end{figure}
      {\ }\\        
    Here are some notes for the table:
    
     {\ }
     
                 \begin{itemize}
               \item [-] The column for $\dim (C)=12$ uses \eqref{nstreq}, while others rely on \eqref{streq}.
               
           {\ }
                
              \item [-]   The asterisk ***  in the table indicates that a solution of the O.D.E exists 
in some interval 
                  $[0,\theta]$, but the vanishing angle based on the method does
not exist in general.

 {\ }
 
          \item [-]     However, since the equality of
                    \eqref{str}\text{ and }\eqref{nstr}
may not be attained,
slots filled with ``***" could still support actual vanishing angles.

{\ }
 
\item [-]    Values in parentheses are tentative, pending on more numerical analysis.
                  \end{itemize}
                     
{\ }

              Generally, using \eqref{streq} gives more accurate data.
              Nevertheless, \eqref{nstreq} has a key advantage,
              Proposition 1.4.2 in \cite{Law}:
 
 {\ }
 
             {\it
              Let $k$ be the dimension of a minimal cone under consideration,
              $\alpha$ explained in \eqref{aofstr},
              and
              $\theta_c(k,\alpha)$ the vanishing angle for \eqref{nstreq} depending on $k\text{ and }\alpha$.
              Then for any $\ell>k$,
              \begin{equation}\label{ind}
              \tan(\theta_c(\ell, \frac{\ell}{k}\alpha))<\frac{k}{\ell}\tan(\theta_c(k,\alpha)).
              \end{equation}
              }
              
               {\ }

              One can use this property to estimate vanishing angles for $\dim (C)>12$.
            Another useful property is the following monotonicity:
                  \begin{equation}\label{mono}
                 0< \theta_c(k, \alpha)<\theta_c(k, \beta)\ \text{ for }\ 0\leq \alpha < \beta
                  \end{equation}
             which comes from the monotonicity of $(1-\alpha t)e^{\alpha t}$ in $\alpha$.
             Similarly, one has $\theta_F(k, \alpha)$ for \eqref{streq} and the monotonicity
              \begin{equation}\label{mono2}
                 0< \theta_F(k, \alpha)<\theta_F(k, \beta)\ \text{ for }\ 0\leq \alpha < \beta.
                  \end{equation}

   {\ }

              Now we are ready to proceed.
              
              {\ }

\subsection{Case for $g=4$}\label{g4}

       Let $M_+$ and $M_-$ 
       be focal submanifolds of an isoparametric foliation with {$g=4$}.
       The main points
{of \textbf{Proof of Theorem 1}}
       are:
       \\{\ }

        1. $\frac{\pi}{4}$-normal wedges of $C(M_+)$ (or $C(M_-)$) are nonintersecting.
         \\{\ }

         2.
         In the isoparametric setting, $\theta_0$ only depends on $k,\alpha$
         and we can check everything merely at one point.
         If vanishing angle $\theta_0$ exists and moreover is less than $\frac{\pi}{4}$,
         then 
{the cone}
is minimizing.
                \\{\ }

                3. 
                            Note that $\frac{\sqrt{x-2}}{x}$ is strictly decreasing when $x\geq 4$.
                            So, for $\ell>k
{\geq4}$, we have
                            \[
                            \frac{\sqrt{\ell-2}}{\ell}<\frac{\sqrt {k-2}}{k}\Rightarrow \sqrt {\ell-2}<\frac{\ell}{k}\sqrt {k-2}.
                            \]
                            By \eqref{ind}, \eqref{mono} and the table of upper bounds of vanishing angles, it follows that
                             \begin{equation}\label{ineqs}
                                          \begin{split}
                                 \tan(\theta_c(\ell,\sqrt {\ell-2}))&<\tan(\theta_c(\ell,\frac{\ell}{k}\sqrt {k-2}))
                                                                       <\frac{k}{\ell}\tan(\theta_c(k,\sqrt {k-2}))\\
                                                                          \specialrule{0em}{3pt}{3pt} 
                                    ({when }\ \ell>k=12)  {\ \ \ \ \ \ \ \  }                                 
                                   &<\tan(\theta_c(12,\sqrt{10}))
                                                                       <\tan(9^\circ)
                                                                       <1.
                                             \end{split}
                                \end{equation}
                                {\ }\\{\ }

            4. Hence, $\theta_0<\theta_c(k,\alpha)<45^\circ$ for $k\geq 12$ and $\alpha^2\leq k-2$.
            According to Lawlor's table,
             $0<\theta_0\leq \theta_F(k,\alpha)<45^\circ$
             for $7\leq k\leq 11$ and $\alpha^2\leq k-2$.
           Thus, when $k\geq 7$ and $\alpha^2\leq k-2$,
           \begin{equation}\label{45}
           0<\theta_0<45^\circ.
           \end{equation}
                              \\{\ }

            5. For $g=4$, $M_+$ is of dimension $m_1+2m_2$, so
                     $
                     k_0=\dim(C(M_+))=m_1+2m_2+1.
                     $
            A very delightful
            property of focal submanifolds is
            that,
            for a point $p$ of $M_+$ and any unit normal $\nu$ {at $p$ to $M_+$ in the sphere},
            the second fundamental form
            for $\nu$ with respect to certain orthonormal basis of  $T_xM_+$
            is {(cf. \cite{CR})}
                         \begin{equation} \label{sff}
                         \left(
                         \begin{array}{ccc}
                         I_{m_2} & 0 & 0\\
                        0 & -I_{m_2} & 0\\
                        0 & 0 & O_{m_1}
                         \end{array}
                         \right)  .
                         \end{equation}
{\ }\\
                 Hence
                 \begin{equation}\label{2}
                  \alpha_0^2=2m_2\ \Longrightarrow\ \alpha_0=\sqrt{2m_2}=\sqrt {k_0-1-m_1}\leq\sqrt {k_0-2}.
                  \end{equation}
                  Therefore, when $k_0\geq 7$,
                  by \eqref{2}, \eqref{ineqs} and \eqref{45},
                  the curvature criterion applies and consequently $C(M_\pm)$ are minimizing.
                \\{\ }

               6. When $(m_1, m_2)=(1,2)$,
               $(k,\alpha^2)$ equals $(6,4)$ for $C(M_+)$
               and $(5,2)$ for $C(M_-)$.
               The vanishing angle exists for the latter and $<27^\circ$,
               whereas the former seems subtle
               {because of}
               encountering ``(***)".
               Using the fact $\inf_\nu\det(1-th_{ij}^\nu)=(1+t)^2(1-t)^2$ instead of the control $F(2, t, 5)$,
               one can figure out that vanishing angle for the former exists and is less than $25^\circ$.

               {\ }

               Thus, our proof of Theorem \ref{thm1} is complete.\hfill$\Box$

               {\ }



\subsection{Cases for $g=3$ and $6$}\label{g36}

In both cases, $m_1=m_2=m$. 
Theorem 2 states precisely the following.
{\ }\\{\ }\\
\textbf{Theorem 2$'$.}
          For $(g,m)=(3,2),(3,4),(3,8)$ or $(6,2)$, cones over focal submanifolds are minimizing.
          \begin{proof}
          With respect to a unit normal vector at a point in a focal submanfold, the second fundamental form
          (shape operator) would be similar to \eqref{sff}.
          They are (cf. \cite{CR})

\begin{equation}\label{g3}
\begin{pmatrix}
\dfrac{1}{\sqrt{3}}\cdot I_2 &0\\
0 & -\dfrac{1}{\sqrt{3}}\cdot I_2
\end{pmatrix}
,
\begin{pmatrix}
\dfrac{1}{\sqrt{3}}\cdot I_4 &0\\
0 & -\dfrac{1}{\sqrt{3}}\cdot I_4
\end{pmatrix}
,
\begin{pmatrix}
\dfrac{1}{\sqrt{3}}\cdot I_8 &0\\
0 & -\dfrac{1}{\sqrt{3}}\cdot I_8
\end{pmatrix}
\end{equation}
and
\begin{equation}\label{g6}
\begin{pmatrix}
\sqrt 3 \cdot I_2 & 0 & 0 & 0 & 0\\
0 & -\sqrt 3 \cdot I_2 & 0 & 0 & 0\\
0 & 0& \dfrac{1}{\sqrt{3}}\cdot I_2 &0 &0\\
0 & 0 & 0& -\dfrac{1}{\sqrt{3}}\cdot I_2 &0\\
0 & 0 & 0& 0& O_2
\end{pmatrix}
\end{equation}
respectively.
Hence $(k,\alpha^2)=(5, \frac{4}{3}),(9, \frac{8}{3}), (17,\frac{16}{3})$ and $(11, \frac{40}{3})$.
According to Lawlor's table and \eqref{ind}, vanishing angle $\theta_0$ exists and $<30^\circ$ for each of them.
Namely, the $\theta_0$-normal wedges over the focal submanifold under consideration are nonintersecting.
So the statement stands.
          \end{proof}
\begin{rem}
For $(g,m)=(3,1)$, $(k,\alpha^2)=(3,\frac{2}{3})$ and 
      $(g,m)=(6,1)$, $(k,\alpha^2)=(6,\frac{20}{3})$,
      we cannot apply the estimate in the table.
\end{rem}          
%

{\ }

\section{Cones over products of focal submanifolds}\label{Section4}

{We introduce the minimal product of two minimal submanifolds in spheres
by illustration for focal submanifolds
 in \S\ref{mp}.
Detailed estimates for the normal radius and the vanishing angles are exhibited in \S\ref{NR} and \S\ref{EVA} respectively.
To apply the curvature criterion we make comparisons in \S\ref{CP}.
Since the minimal product can also be defined for multiple minimal submanifolds,
considerations are given for that situation in \S\ref{added}.
Finally, the case involving $g=2$, for which focal submanifolds are geodesic spheres, 
is discussed in \S\ref{g2}.
}
\subsection{Minimal product}\label{mp}
         Given two focal submanifolds $f_1:M_1^{k_1}\hookrightarrow S^{n_1}$ and $f_2:M_2^{k_2}\hookrightarrow S^{n_2}$ for isoparametric foliations of unit spheres, with $g_1,\, g_2$ respectively,
         define an embedding $G: M\triangleq M_1^{k_1}\times M_2^{k_2} \rightarrow S^{n_1+n_2+1}$ by
         \[
         (x,y)\mapsto (\lambda f_1(x), \mu f_2(y)), \text{  \    \       with\ \ } \lambda=\sqrt{\dfrac{k_1}{{k_1+k_2}}}\ \text{\ \ and }\ \mu=\sqrt{\dfrac{k_2}{{k_1+k_2}}}.
         \]
         We write $x$ and $y$ short for $f_1(x)$ and $f_2(y)$.
         Together with
         $
         \eta_0=(\mu x, -\lambda y)
        $,
         orthonormal bases
         $\{\sigma_1,\cdots, \sigma_{n_1-k_1}\}$ and $\{\tau_1,\cdots,\tau_{n_2-k_2}\}$
         of $T_x^\perp M_1$ and $T_y^\perp M_2$ 
         induce an orthonormal basis
         $\{(\sigma_1,0)$,$\cdots, (\sigma_{n_1-k_1},0)$,
         $(0,\tau_1)$,$\cdots$, $(0,\tau_{n_2-k_2})$, $\eta_0\}$
         of the normal space of $G(M)$ at $P\triangleq(\lambda x,\mu y)$ in $S^{n_1+n_2+1}$.

         Let $A$ be the symbol of shape operators.
         Then $A_{(\sigma_i,0)},\; A_{(0,\tau_j)}:T_PM\rightarrow T_PM$,\, $A_{\sigma_i}:T_xM_1\rightarrow T_xM_1$
         and $A_{\tau_j}:T_yM_2\rightarrow T_yM_2$
         have the following relations:

         \[ A_{(\sigma_i,0)} =
\begin{pmatrix}
\dfrac{1}{\lambda}A_{\sigma_i} & O\\
O &O
\end{pmatrix},
\]
         and
         \[ A_{(0,\tau_j)} =
\begin{pmatrix}
O & O\\
O & \dfrac{1}{\mu}A_{\tau_j}
\end{pmatrix}.
\]
{\ }\\
          Also note that there are three Levi-Civita connections
            \begin{eqnarray}
                     M\ \longhookrightarrow &S^{n_1+n_2+1}\longhookrightarrow &\mathbb R^{n_1+n_2+2}
                       \nonumber\\
                          \specialrule{0em}{2pt}{2pt} 
                    \nabla\ \ \ \ \ \ \ &\overline{\nabla} &\ \ \ \ \ \ \ \mathcal D
                       \nonumber
                       \end{eqnarray}
         and, for tangent vector fields $X=(X_1,X_2)$ and $Y=(Y_1,Y_2)$ around $P$, we have
              \begin{eqnarray*}
                     <A_{\eta_0} X, Y>
                     &=&
                     <\overline{\nabla}_XY-{\nabla}_XY,\ \eta_0>\\
                   &=& <\mathcal D_XY,\ \eta_0>\\
                    &=& -<\mathcal D_X\eta_0,\ Y>\\
                      &=&-<\mathcal D_{(X_1,X_2)}(\mu x
{,}
                      -\lambda y),\ Y>\\
                      &=& <(-\frac{\mu}{\lambda}X_1,\frac{\lambda}{\mu}X_2),\; (Y_1,Y_2)>.
                       \end{eqnarray*}
         Therefore, $$A_{\eta_0}=
                 \begin{pmatrix}
                 -\dfrac{\mu}{\lambda}I_{k_1} & O\\
                 O & \dfrac{\lambda}{\mu}I_{k_2}
                 \end{pmatrix}
                 $$
         and trace$(A_{\eta_0})=0$.
         These imply that $G$ minimally embeds $M$ into $S^{n_1+n_2+1}$.
{We call $G(M)$ the minimal product of $M_1$ and $M_2$.}
\\{\ }

\subsection{For normal radius}\label{NR}
       At $P=(\lambda x, \mu y)\in M$,
         let $$N=\big(a_1\xi_1 + a_0 \mu x, a_2\xi_2 - a_0\lambda y\big),
         \, \,\,\, a_0,\, a_1,\, a_2\in \mathbb R$$
         where $\xi_1$ is a unit normal to $M_1$ at $x$ in $S^{n_1}$
         and $\xi_2$ a unit normal to $M_2$ at $y$ in $S^{n_2}$,
         such that
                $\omega\triangleq \|N\|$ attains the smallest (nonzero) for
                \[
                Q=P+N\; \in \text{ the cone }C(M).
                \]
  \begin{figure}[ht]
 \includegraphics[scale=0.75]{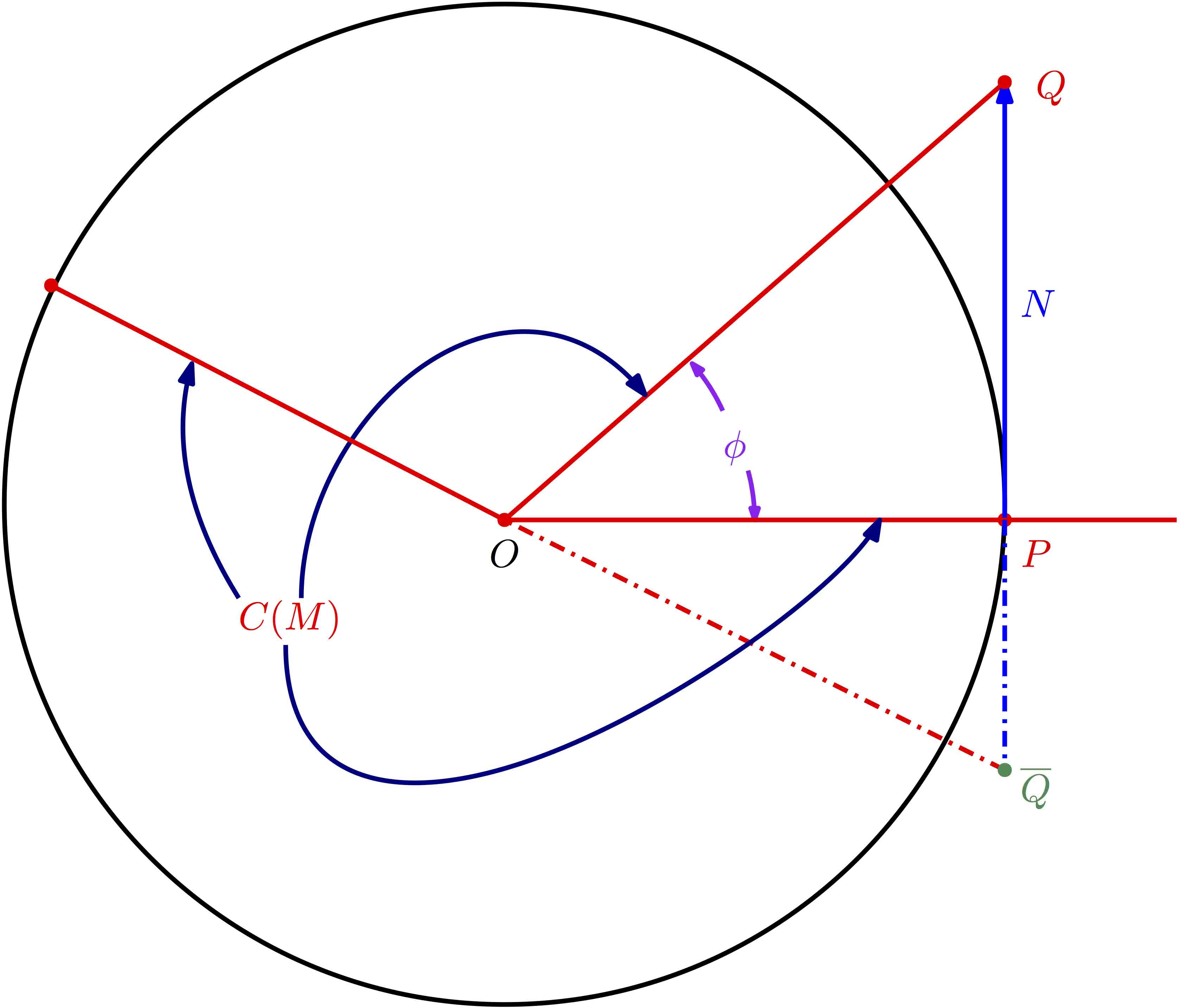}
\caption{{$Q$ and fake $\overline{Q}$}}
 \end{figure}
     {\ }
     
     Let $\phi\triangleq\arctan \omega$.
There are three possibilities to consider for smallest $\omega$.
If none of
{
\textbf{(I),\ (II)} and \textbf{(III)}
below
give a finite $\omega$}, we define $\phi=\frac{\pi}{2}$.
               Suppose $\phi<\frac{\pi}{2}$
               (the actual situation for $g_1,\ g_2\in\{3,\ 4,\ 6\}$).
               Then $Q=\Big(a_1\xi_1 + (a_0 \mu+\lambda) x,\; a_2\xi_2+( - a_0\lambda +\mu)y\Big)$.

{\ }

                \textbf{(I)}. $a_1\cdot a_2\neq 0$.
                
                Recall that (\cite{CR}) along any normal geodesic of $M_1$ in $S^{n_1}$,
                it takes length $\frac{\pi}{g_1}$ to arrive at its pairing focal submanifold.
                So before $\frac{2\pi}{g_1}$ it cannot get back to $M_1$.
                The same holds for $M_2$ in $S^{n_2}$.
{As a consequence,}
                  \begin{equation*}
                               \begin{cases}
                                        \  \dfrac{a_1^2+(a_0\mu+\lambda)^2}{\lambda^2}\ =\ \dfrac{a_2^2+(\mu-a_0\lambda)^2}{\mu^2};\ \ \ & \\
                                      \    \ \ \ \ \   \ \ \ \ \ \dfrac{a_0\mu+\lambda}{a_1}\ =\ \cot \phi_1 ,&\ \  \text{ where }\phi_1\in \cfrac{2\pi}{g_1}\cdot \mathbb Z\ ;
                                          \\
                                        \  \ \ \ \ \  \ \ \ \ \ \dfrac{\mu-a_0\lambda}{a_2}\ =\ \cot \phi_2,&\ \ \text{ where }\phi_2\in \cfrac{2\pi}{g_2}\cdot \mathbb Z\ .
                                          \\
                                           \end{cases}
                       \end{equation*}
                       Since $a_1a_2\neq 0$,
                       it indicates that $\phi_1,\, \phi_2$ cannot belong to $\pi\cdot \mathbb Z$.
                       Therefore, we can have
                       $$\cfrac{a_1^2}{\lambda^2\sin^2\phi_1}=\dfrac{a_2^2}{\mu^2\sin^2\phi_2}\triangleq t^2,$$
                       so that
                       $ a_1=t{\lambda}\sin \phi_1$ and $a_2=\pm t{\mu}\sin \phi_2$.
                       Thus,
                \begin{equation*}
                               \begin{cases}
                                         a_0\mu + \lambda =\ \ \ t\lambda\cos \phi_1\\
                                         \mu -a_0 \lambda = \pm t\mu\cos \phi_2\\
                                           \end{cases}
                                           \end{equation*}
                                          and consequently, by the relation $\lambda^2+\mu^2=1$, we have
                              \begin{equation}\label{4.1}
                                             \begin{cases}
                                       {\; \:}  1= t(\lambda^2\cos \phi_1\pm \mu^2\cos\phi_2),\\
                                         a_0=t\lambda\mu(\cos\phi_1\mp \cos\phi_2). \ \\
                                           \end{cases}
                       \end{equation}
                       Now
                       \begin{eqnarray}
                       \omega^2 
                       =
                       a_1^2+a_2^2+a_0^2
                       =
                       t^2[\lambda^2\sin^2\phi_1+\mu^2\sin^2\phi_2+\lambda^2\mu^2(\cos\phi_1\mp \cos\phi_2)^2],
                       \nonumber
                       \end{eqnarray}
                                    which, together with \eqref{4.1}, produces
                        \begin{equation}\label{4.2}
                       \omega^2 (\lambda^2\cos \phi_1\pm \mu^2\cos\phi_2)^2
                       = [\lambda^2\sin^2\phi_1+\mu^2\sin^2\phi_2+\lambda^2\mu^2(\cos\phi_1\mp \cos\phi_2)^2].
                         \end{equation}
         By adding $(\lambda^2\cos \phi_1\pm \mu^2\cos\phi_2)^2$ to both sides of \eqref{4.2}, we have
                         \begin{eqnarray}\label{4.3}
                       & &(\omega^2+1) (\lambda^2\cos \phi_1\pm \mu^2\cos\phi_2)^2\nonumber\\
                       &=&[\lambda^2\sin^2\phi_1+\mu^2\sin^2\phi_2+\lambda^2\cos^2\phi_1+\mu^2\cos^2\phi_2]\\
                       &=&1.\nonumber
                       \end{eqnarray}
                To seek for the smallest $\omega$, we consider the chance
                  \begin{equation}\label{phiA}
                  \omega^2=
                  \left(
                  \dfrac{1}{\lambda^2|\cos\frac{2\pi}{g_1}|+\mu^2|\cos\frac{2\pi}{g_2}|}
                  \right)^2-1.
                       \end{equation}
                       Notice that $g_1,\, g_2\in \{3,\ 4,\ 6\}$.
                       By ranges of $\phi_1$ and $\phi_2$,
                   \begin{equation}\label{A}
                  \omega^2
                  \geq
                  \left(
                  \dfrac{1}{\lambda^2\cdot \frac{1}{2}+\mu^2\cdot \frac{1}{2}}
                  \right)^2-1=3.
                       \end{equation}
                       
                 {\it {Remark}. }In fact, when $g_1=g_2=4$, \eqref{phiA} generates infinity.
                 However, in this case,
                 the realizable \eqref{phiB} below asserts the finiteness of $\tan^2\phi$.

  {\ }

                       \textbf{(II)}. $a_1\neq0,\; a_2=0$. (Similar for $a_1=0,\; a_2\neq0$.)

We have {two restrictions}
                        \begin{equation*}
                               \begin{cases}
                                          \dfrac{
                                       \   a_1^2+(a_0\mu+\lambda)^2}{\lambda^2}\ =\ \dfrac{(\mu-a_0\lambda)^2}{\mu^2},\ \ \ &
                                                                                \\
                                                                                   \specialrule{0em}{2pt}{2pt} 
                                       \   \ \ \ \ \ \ \ \   \ \ \dfrac{a_0\mu+\lambda}{a_1}\ =\ \cot \phi_1,&
                                      \    \text{ where }\phi_1\in \cfrac{2\pi}{g_1}\cdot \mathbb Z\ .
                                          \\
                                           \end{cases}
                                           \end{equation*}

                      As argued in \textbf{(I)}, $\phi_2$ cannot belong to $\pi\cdot \mathbb Z$ and assume
                        \begin{equation}
                        \dfrac{a_1^2}{\lambda^2\sin^2\phi_1}
                        =
                        \cfrac{(\mu-a_0\lambda)^2}{\mu^2}
                        \triangleq
                        t^2,
                        \end{equation}
                        so that
                  \begin{equation}
                               \begin{cases}
                                                                         a_1=\ \ t\lambda\sin\phi_1,\\
                                         \mu -a_0 \lambda =\ \pm  t\mu.
                                           \end{cases}
                   \end{equation}
                      Similarly,
                       \begin{equation}
                                                                a_0\mu + \lambda = t\lambda\cos \phi_1
                                                                    \end{equation}
                 and
                 \begin{equation}\label{4.9}
                               (\omega^2+1) (\lambda^2\cos\phi_1\pm \mu^2)^2=1.
                 \end{equation}
                            Again, for the smallest $\omega$, we focus on the likelihood
                              \begin{equation}\label{phiB}
                  \omega^2=
                  \left(
                  \dfrac{1}{\lambda^2|\cos\frac{2\pi}{g_1}|+\mu^2}
                  \right)^2-1.
                       \end{equation}
By the range of $\phi_1$,
                                                                                  \begin{equation}\label{B}
                  \omega^2
                  \geq
                  \left(
                  \dfrac{1}{\lambda^2\cdot \frac{1}{2}+\mu^2}
                  \right)^2-1.
                       \end{equation}                                          

                       {\ }

                     \textbf{(III)}. $a_1=a_2=0$.

                     In this case,
                     $Q=\Big((a_0 \mu+\lambda) x,\; ( - a_0\lambda +\mu)y\Big)\in C(M)$
                     which implies
                      \begin{equation}\label{4.10}
                      \dfrac{(a_0\mu+\lambda)^2}{\lambda^2}= \dfrac{(\mu-a_0\lambda)^2}{\mu^2}.
                       \end{equation}
  Recall that
 $\lambda=\sqrt{\dfrac{k_1}{{k_1+k_2}}}$\  and\  $\mu=\sqrt{\dfrac{k_2}{{k_1+k_2}}}.$
            It is not hard to see that,
            only when $\lambda\neq \mu$,
            there is
a finite solution
 \begin{equation*}
                       a_0= \dfrac{2\sqrt{k_1k_2}}{k_1-k_2},
  \end{equation*}
                      and
{possibly}
                                \begin{equation}\label{C}
                      \omega^2=
                      \dfrac{4k_1k_2}{(k_1-k_2)^2}.
   \end{equation}
{\ }

             Hence,
             taking the possible occurrence of fake $\overline Q$
             in Figure 1 into account,
             $\tan^2 \phi$ is no less than,
             by definition,
              the smallest quantity in \eqref{A}, \eqref{B} and \eqref{C}.

  {\ }

 \subsection{Existence of vanishing angles}\label{EVA}
                     Let $\xi_1,\; \xi_2$ be unit normals to $M_1$ and $M_2$ at $x$ and $y$ in $S^{n_1}$ and $S^{n_2}$,
                     and  $A_{\xi_1}$, $A_{\xi_2}$ the corresponding shape operators
                     respectively.
                     Set
                            \begin{eqnarray}
                       \eta_1&=&\big(\xi_1,\;  0\big),\nonumber\\
                       \eta_2&=&\big( 0,\;  \xi_2\big),\nonumber\\
                     \eta_0&=&\big(\mu x,\; -\lambda y\big).\nonumber
                       \end{eqnarray}
                       Then, for a unit normal 
              $\varsigma\triangleq\sum\limits_{i=0}^2 c_i\eta_i$ with $\sum\limits_{i=0}^2 c_i^2=1$,
                       the shape operator 
                $A_\varsigma$
                       of $M$ at $(\lambda x, \mu y)$ in $S^{n_1+n_2+1}(1)$ is 
                       \[
                       \begin{pmatrix}
                      \dfrac{1}{\lambda}(-c_0\mu I_{k_1}+ c_1 A_{\xi_1})& O\\
                       O &  \dfrac{1}{\mu} (c_0\lambda I_{k_2}+c_2 A_{\xi_2})
                       \end{pmatrix}.
                       \]
                       Since $M_1$ and $M_2$ are minimal,
                       trace$(A_{\xi_1})=0$
                       and  trace$(A_{\xi_2})=0$.
                        Set $S\triangleq \dim(M)=k_1+k_2$.
                        It follows that
                        \begin{eqnarray}\label{normctrl}
                        \|A_{\varsigma}\|^2
                        &=&
                        c_0^2 S+
                        c_1^2
                        \dfrac{S}
                        {k_1}\|A_{\xi_1}\|^2+
                         c_2^2
                         \dfrac{S}
                         {k_2}\|A_{\xi_2}\|^2.
                       \end{eqnarray}
                       Observations from \eqref{sff}, \eqref{g3} and \eqref{g6} show that (same for $(m_1,m_2)=(1,1)$)
                             \begin{eqnarray}
                       \text{if } M_i \text{ corresponds to } g=3,\, &
                       \text{then } \|A_{\xi_i}\|^2= \dfrac{k_i}{3};
                       \nonumber\\
                    \text{if } M_i \text{ corresponds to } g=4,\, &
                       \text{  \ \  \ \ then } \|A_{\xi_i}\|^2\leq {k_i-1};
                       \nonumber\\
                     \text{if } M_i \text{ corresponds to } g=6,\, &
                       \text{\; then } \|A_{\xi_i}\|^2= 
                       \dfrac{4}{3}k_i.
                     \nonumber
                       \end{eqnarray}
                       Thus, by  expression \eqref{normctrl},
                     we always have 
            \begin{equation}\label{nsup}
            \|A_{\varsigma}\|^2\leq \frac{4}{3}S.
            \end{equation}
                     
                     {\ }
                     
                     For $S=8,\ 9$ or $10$
                     and $\alpha^2\leq \frac{4}{3}S$,
                     vanishing angle $\theta_F(S+1,\alpha)$ exists according to Lawlor's table
                     and the fact that $\theta_F(9,\sqrt{\frac{32}{3}})
                     <18^\circ.
                     $

           {\ }
                     
                     When $S\geq 11$, we seek for the existence of vanishing angle $\theta_c$.
                   Note that, when the dimension of a cone equals $12$,
                     Lawlor's table confirms the existence of $\theta_c$ for $\alpha^2\leq 19$.
                                          Recall \eqref{ind}: for $\ell>k$
                     \begin{equation}\label{indrepeat}
                       \tan(\theta_c(\ell, \frac{\ell}{k}\alpha))<\frac{k}{\ell}\tan(\theta_c(k,\alpha)).
                     \end{equation}
                     Since $\dfrac{x^2}{x-1}$ is strictly increasing for $x\geq 2$,
                     one can find that, when $\ell>k\geq 2$,
                                \begin{eqnarray}\label{BiLi}
                       \dfrac{(\frac{\ell \alpha}{k})^2}{\ell-1}
                       =
                       \frac{\ell^2 \alpha^2}{k^2(\ell-1)}
            >\dfrac{\alpha^2}{k-1}.
                       \end{eqnarray}
                       Take $(k, \alpha^2)=(12,19)$ and fix $\ell>12$.
                            Then $\theta_c(12, \sqrt{19})$
                       in the latter term of \eqref{indrepeat} exists,
                       and so does $\theta_c(\ell, \frac{\ell}{12}\sqrt{19})$.
                       By monotonicity \eqref{mono},
                       if $\tilde \alpha^2\leq (\frac{\ell \sqrt{19}}{12})^2$,
                       $\theta_c(\ell,\tilde \alpha)$ exists.
                       In particular,
                       according to \eqref{BiLi},
                       the existence of $\theta_c(\ell, \tilde\alpha)$ for
                        $\tilde \alpha^2=\frac{4}{3}(\ell-1)
                        <
                        \frac{19}{11}(\ell-1)<(\frac{\ell \sqrt{19}}{12})^2$ is guaranteed.
                       As a result,
                       $\theta_c\left(S+1, \sqrt{\frac{4}{3}S}\right)$ exists for
 $S\geq 11$.
                                          
      {\ }

Hence, $\theta_F$ and $\theta_0$ always exist for $S\geq 8$.

{\ }

\subsection{Comparison between $\phi$ and $2\theta_0$}\label{CP}
                     If we can get 
                     $$(\star)\ \ \
                      {2}\theta_0<\phi
                      $$                 
{by the relation given at the end of \S\ref{NR},
                     then the curvature criterion applies
                     and the corresponding cone is area-minimizing.
                     
                     For our purpose, let us make some observations about $\phi$.                     
                     By Appendix, {it follows that}
                   \begin{equation}\label{compare}
                    \min\left\{3,\ \dfrac{4k_1k_2}{(k_1-k_2)^2}\,\right\}
                   >
                   \min
                   \left\{
                  \left(
                  \dfrac{1}{\lambda^2\cdot \frac{1}{2}+\mu^2}
                  \right)^2-1
,\ \left(
                   \dfrac{1}{\lambda^2+\mu^2\cdot \frac{1}{2}}
                                     \right)^2-1
                   \right\}.
                   \end{equation}
                   Assume that $k_1\leq k_2$ from now on.
                   One can see that
                   \begin{equation*}
{  \left(
                   \dfrac{1}{\lambda^2+\mu^2\cdot \frac{1}{2}}
                                     \right)^2-1
                   \geq 
                                       \left(
                  \dfrac{1}{\lambda^2\cdot \frac{1}{2}+\mu^2}
                  \right)^2-1.
 }
                   \end{equation*}
                   Hence, 
                   $$\tan^2\phi
                   \geq
{                                                        \left(
                  \dfrac{1}{\lambda^2\cdot \frac{1}{2}+\mu^2}
                  \right)^2
                  }
                  -1,
                                     $$
                   and therefore
                   \begin{equation}\label{nice}
                   \cos\phi\leq
                   1-\dfrac{k_1}{2S}.
                   \end{equation}
    
                   
                 \begin{rem}\label{rknice}
Since $t$ in \S\ref{NR} is allowed to be negative,
our computations include the possibility of fake $\overline Q$.
                   It implies that the upper bound $1-\frac{k_1}{2S}$
                   is in fact uniform for $\left|\cos\phi_\iota\right|$
                   where $\phi_\iota$ corresponds to
                   any intersection point (if existed
                   {other than} $\pm P$)
                   of               
{the submanifold}
                   $G(M_1\times M_2)$ and 
                   a great circle perpendicular to 
                   it through $P$.
                   \end{rem}


Now we can prove a weak version of Theorem 3.
\\{\ }\\
\textbf{Theorem 3$'$.}
       {\it  Cones over the minimal products of {\tt two}
       focal submanifolds of isoparametric foliations
       with {$g=3,4,6$} and $(m_1,m_2)\neq (1,1)$ 
       are area-minimizing.}
{\ }\\{\ }\\
\textbf{Proof.}
                 Note that $4\leq k_1\leq k_2$ in the case.
                 Since $\theta_F(S+1,\sqrt{\frac{4}{3}S})$ exists for $S\geq 8$,
we have, by \eqref{str}, \eqref{Fofstr}, \eqref{streq}, \eqref{mono2} and \eqref{nsup}, that
\begin{equation}\label{4c}
                  \tan \theta_0
                  \leq
                  \tan \theta_F(S+1,\sqrt{\frac{4}{3}S})
                  <
                  \frac{1}{\sqrt{\frac{4}{3}S}}\sqrt{\dfrac{S}{S-1}}=\dfrac{\sqrt 3}{2}\sqrt{\dfrac{1}{S-1}}.
                  \end{equation}
                  Then
                   \begin{equation*}
                  \tan 2\theta_0< \dfrac{4\sqrt{3(S-1)}}{4S-7}.
                  \end{equation*}
                  According to
                  \begin{equation}\label{thm3tanphi}
                  \begin{split}
                  \tan^2\phi
                 & \geq
                   \left(\dfrac{1}{1-\dfrac{k_1}{2S}}\right)^2-1
                  \geq \left(\dfrac{1}{1-\dfrac{4}{2S}}\right)^2-1\\
                  \specialrule{0em}{5pt}{5pt} 
              &=\dfrac{4(S-1)}{(S-2)^2}>\dfrac{3(S-1)}{(S-\frac{7}{4})^2}
                  >
                  \tan^2 2\theta_0\ ,
                  \end{split}
                  \end{equation}
we have $(\star)$ and thus complete the proof.
 \hfill$\Box$
\hspace{0.2cm}
{\ }\\{\ }

                   Let us assume $S\geq11$ and figure out when $(\star)$ holds.
                   Note that 
                   $$\tan\left[\theta_c(12,\sqrt{\frac{44}{3}})\right]<\tan 9.55^\circ<0.1683
                   .$$
                   Combined with \eqref{indrepeat}, \eqref{BiLi} and \eqref{mono}, 
                   the existence of $\theta_c\Big(S+1,\ \sqrt{\dfrac{4}{3}S}\Big)$ leads to
                   \begin{equation}\label{new}
                   \tan\theta_0
                   <
                   \tan\left[
                  \theta_c\Big(S+1,\ \sqrt{\dfrac{4}{3}S} \Big)\right]
                   <
                   \dfrac{12}{S+1}\tan\left[
                  \theta_c\Big(12,\ \sqrt{\dfrac{44}{3}} \Big)\right].
                   \end{equation}
                   Then
                         \begin{equation*}
                          \tan2\theta_0
                          <
                          \dfrac{2\cdot \dfrac{ 12 }{S+1}\cdot 0.1683}{1-\left(\dfrac{ 12 }{S+1}\cdot 0.1683\right)^2}\ .
                 \end{equation*}
                     Hence, with the agreement that $k_1\leq k_2$, it is sufficient to {see when}
                     \begin{equation}\label{thm4tanphi}
                     \tan^2\phi
                     \geq
                      \left(\dfrac{1}{1-\dfrac{k_1}{2S}}\right)^2-1
                     \geq
                     \left(\dfrac{2\cdot \dfrac{ 12 }{S+1}\cdot 0.1683}{1-\left(\dfrac{ 12 }{S+1}\cdot 0.1683\right)^2}\right)^2
                     >
                     \tan^2 2\theta_0\ ,
                     \end{equation}
{or equivalently,}
                    \begin{equation}\label{poly}
                    \begin{split}
                    & \left(k_1S-\frac{k_1^2}{4}\right)
                     \left[(S+1)^2
                     -(12*0.1683)^2
                     \right]^2\\
                    & {\ \ \ \ \ \ \ \ \ \ \ \ \ \ \ \ \ \ \ \ \ \ \ \ \ \ \ \ \    \ \ \ \ } -
                     \left(
                     S-\frac{k_1}{2}
                     \right)^2
                     \left[
                     24(S+1)*0.1683
                     \right]^2
{\geq}
                     0.
                     \end{split}
                     \end{equation}
       
       {\ }              

Based on the analysis, we get a weak version of Theorem 4.\\{\ }\\
  \textbf{Theorem 4$'$.}
       {\it  Cones, of dimension no less than $10$, over the minimal products
       of 
       {\tt two} 
       focal submanifolds
       with {$g=3,4,6$}
       are area-minimizing.}
{\ }\\{\ }\\ 
\textbf{Proof.}
             When $k_1\geq 4$, a proof is given in that of Theorem 3.
             Consider now that $k_1=2$ and $3$.
             To determine the sign of polynomial \eqref{poly}, we run the following codes in Mathematica: 
                         \begin{figure}[h]
                        \begin{center}
                        \includegraphics[scale=0.29]{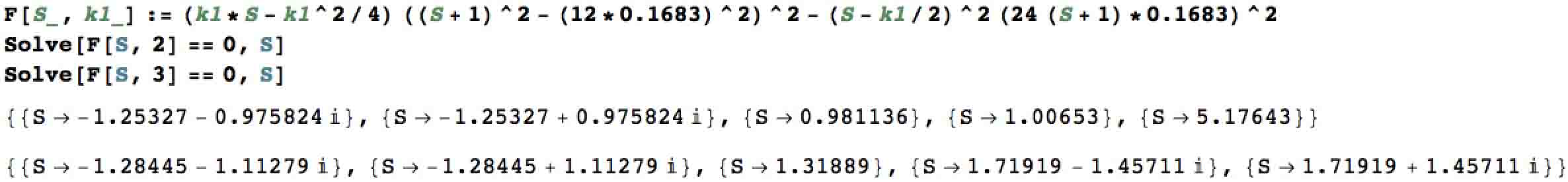}
                        \end{center}
              \end{figure}
              \\
                   and it can be seen that, whenever $S\geq 11$ for \eqref{new}, we have \eqref{poly} hold.
              
              Let us check the cases of $S=9,\ 10$ for $k_1=2$ (which can imply the cases for $k_1=3$).
              By Lawlor's table,
              $\theta_F(10,\sqrt{12})<13.51^\circ$
              and
              $\theta_F(11,\sqrt{\frac{40}{3}})<11.35^\circ$.
              Then direct computations show that
              $\tan^2\phi\geq (\frac{1}{1-\frac{1}{9}})^2-1>\tan^2 27.02^\circ$
              and 
              $\tan^2\phi\geq (\frac{1}{1-\frac{1}{10}})^2-1>\tan^2 22.70^\circ$ respectively.
              So, $\phi>2\theta_F$ and the proof gets complete.
               \hfill$\Box$
\hspace{0.2cm}

{\ }

              %

\subsection{Extension to the multiple case}\label{added}
          In this subsection, we explain ideas
          of extending Theorem 3$'$ and Theorem 4$'$
          to products of three focal submanifolds with restriction of each $g\in\{3,\, 4,\, 6\}$.
          Case involving more components
          can be similarly done by induction.
          In particular, Theorem 3 
           will be proved.
          
           Given three focal submanifolds
           $f_1:M_1^{k_1}\hookrightarrow S^{n_1}$,
           $f_2:M_2^{k_2}\hookrightarrow S^{n_2}$ and 
           $f_3:M_3^{k_3}\hookrightarrow S^{n_3}$
           for isoparametric foliations of unit spheres.
          %
         Define $\tilde G: \tilde M\triangleq
         M_1^{k_1}\times M_2^{k_2}\times M_3^{k_3}
         \rightarrow S^{n_1+n_2+n_3+2}(1)$ by
         \[
         (x,y,z)\mapsto (\lambda' f_1(x),\ \mu' f_2(y),\ \tau' f_3(z)), \text{  \           with }\]
         \[ \lambda'=\sqrt{\dfrac{k_1}{{k_1+k_2+k_3}}},\ \  \
         \mu'=\sqrt{\dfrac{k_2}{{k_1+k_2+k_3}}}\  \  \text{ and }\ 
          \tau'=\sqrt{\dfrac{k_3}{{k_1+k_2+k_3}}}.
         \]
         Then $\tilde G: \tilde M \hookrightarrow S^{n_1+n_2+n_3+2}(1)$ is a minimal embedding,
         which follows from that
         \[
         M_1\times M_2\times M_3\cong (M_1\times M_2)\times M_3,
         \text{\ \ \  and \ \ } \tilde G(x,\ y,\ z)=G(G(x,\ y), \ z)
         \]
         where map $G$ is given for two minimal submanifolds at the beginning of this section.
         
         Let $\tilde S=k_1+k_2+k_3$.
         With the above understanding,
         it easily follows, from expressions \eqref{nsup} and \eqref{normctrl},
         that
         the norm square of shape operator of $\tilde M$
         at every point for any unit normal is no more than $\frac{4}{3}\tilde S$.
          Hence, vanishing angles $\tilde \theta_F$ and $\tilde \theta_0$
          corresponding to that in \S\ref{EVA}
          and those 
          in proofs of Theorem 3$'$ and Theorem 4$'$
          exist.
          
          Without loss of generality, assume that $k_1\leq k_2\leq k_3$.
          Let 
          $\phi$ and $\tilde \phi$ be the normal radii of
          $M_1\times M_2$ in $S^{n_1+n_2+1}(1)$ 
          and $\tilde M$ in $S^{n_1+n_2+n_3+2}(1)$ respectively.
          We stick to the notation $S=k_1+k_2$.
          Set
                  $\tilde \lambda=\sqrt{\frac{S}{\tilde S}}$
 {and}
                  $\tilde \mu=\sqrt{\frac{k_3}{\tilde S}}$.
          Then, applying the discussions of \textbf{(I)},\ \textbf{(II)}\ and \textbf{(III)}
          to $(M_1\times M_2)\times M_3$,
          we have, according to Appendix and Remark \ref{rknice} (in \S \ref{CP}), that
\begin{equation}\label{multi}
\begin{split}
     \tan^2\tilde\phi
          & \geq 
 {
   \min
           \left\{
                      \left(
                        \dfrac{1}{\tilde \lambda^2
                        \cos \phi
                        +\tilde \mu^2}
                        \right)^2-1,\
                              \left(
                  \dfrac{1}{\tilde \lambda^2
                       +\tilde \mu^2\cdot\frac{1}{2}}
                        \right)^2-1
            \right\}
 }\\
    &     \geq
         \left(
                    \dfrac{1}{\tilde \lambda^2
                      (1-\frac{k_1}{2S})
                +\tilde \mu^2}
                     \right)^2-1.
\end{split}
\end{equation}
      Hence,
      we gain an inheritable relation
                           \begin{equation}\label{beautiful}
                      \cos\tilde\phi
                      \leq
                      1-\frac{k_1}{2\tilde S}
                      \end{equation}
where the {upper bound has the uniform} property
{as} in Remark \ref{rknice}.
             Based on the proofs of Theorem 3$'$ and Theorem 4$'$,
              it then follows correspondingly that
               \begin{equation}\label{Qnice}
                      \tan^2\tilde\phi
                      \geq
                      \left(\dfrac{1}{1-\dfrac{k_1}{2\tilde S}}\right)^2-1
                      >\tan^2 2\tilde \theta_F\geq \tan^2 2\tilde \theta_0.
                      \end{equation}

           The descendent \eqref{beautiful} and
{the} coupled property in Remark \ref{rknice} 
           are crucial for this procedure to be valid
           for the induction
           on the number of focal submanifolds.
   
    {\ }
           
\subsection{To include $g=2$}\label{g2}
         Due to the distinct behaviors of corresponding focal submanifolds,
         we save the discussions about $g=2$ separately
         and accomplish the proof of Theorem 4
          in this subsection.
          Since the case of products
{purely of} spheres has been classified by \cite{Law},
          we only consider the mixed type.
          More explicitly,
          let $M_1^{k_1},\ \cdots,\ M_r^{k_r}$ be focal submanifolds in spheres corresponding to $g\in\{3,\ 4,\ 6\}$,
          and $N_1^{l_1},\ \cdots,\ N_s^{l_s}$ focal submanifolds for $g=2$.
          Then $M=M_1^{k_1}\times \cdots \times M_r^{k_r}\times N_1^{l_1}\times\cdots\times N_s^{l_s}$
          with $r,\ s \geq 1$ is what we are concerned about.
          
        {\ }

          \textbf{(A).\ }
          $M=(M_1^{k_1}\times \cdots \times M_r^{k_r})\times N_1^{l_1}$. Namely $r\geq 1$ and $s=1$.
          Denote by $\hat \phi$ the normal radius of the minimal submanifold $M$ in the sphere.
          Assume that $k_1\leq\cdots\leq k_r$ and
          set $K=k_1+\cdots+k_r$.
          Then \textbf{(I)}, \textbf{(II)} and \textbf{(III)}
          give a lower bound 
          $$
          \min
           \left\{
                      \left(
                        \dfrac{1}
                        {
                        \frac{K}{K+l_1}
                        (
                        1-\frac{k_1}{2K}
                        )
                        +\frac{l_1}{K+l_1}
                        }
                        \right)^2-1,\
         \dfrac{4Kl_1}{(K-l_1)^2}
            \right\}.
            $$
          of $\tan^2\hat\phi$.  
          Although the last term of the following
          has no apparent geometric meaning in the current case,
          one can deduce from Remark \ref{stronger} that 
                   $$
                   \tan^2\hat\phi
                   \geq
                   \min
           \left\{
                      \left(
                        \dfrac{1}
                        {
                       1-\frac{k_1}{2(K+l_1)}
                        }
                        \right)^2-1,\
                              \left(
                        \dfrac{1}
                        {
                       1-\frac{l_1+1}{K+l_1}
                        }
                        \right)^2-1
            \right\}.
            $$
Therefore, 
                             \begin{equation}\label{g2s1}
                   \tan^2\hat\phi
                   \geq
                      \left(
                        \dfrac{1}
                        {
                       1-\dfrac{\min\left\{k_1, 2(l_1+1)\right\}}{2(K+l_1)}
                        }
                        \right)^2-1.
            \end{equation}
 {\ }
            
            \textbf{(B).\ }
            $M=(M_1^{k_1}\times \cdots \times M_r^{k_r})\times (N_1^{l_1}\times\cdots\times N_s^{l_s})$ with $r\geq 1$ and $s\geq 2$.
            We use symbols $\phi_0$ and  $\hat \phi$
            to represent the normal radii of minimal products
            $N_1^{l_1}\times\cdots\times N_s^{l_s}$ and $M$
            in the corresponding spheres, respectively.
            Assume that $l_1\leq \cdots \leq  l_s$ and set $L=l_1+\cdots+l_s$.
            Then, according to \S 5.1 of \cite{Law},
            $\phi_0=\cos^{-1}\left(1-\frac{2l_1}{L}\right)$.
            Also note that
            the minimal
            $N_1^{l_1}\times\cdots\times N_s^{l_s}$
             is symmetric about the origin.
             Hence,
             $\cos\phi_0=1-\frac{2l_1}{L}$ has the uniform upper bound property mentioned in Remark \ref{rknice}.
             Therefore,
                \begin{equation} \label{g2s2}
         \begin{split}
                                       &{\ \ \ }\tan^2\hat\phi\\
                   &\geq
                        \min
           \left\{
                      \left(
                        \dfrac{1}
                        {
                        \frac{K}{K+L}
                        (
                        1-\frac{k_1}{2K}
                        )
                        +\frac{L}{K+L}
                        }
                        \right)^2-1,\
         \left(
                        \dfrac{1}
                        {
                        \frac{K}{K+L}
                        +\frac{L}{K+L}
                         (
                        1-\frac{2l_1}{L}
                        )
                        }
                        \right)^2-1
            \right\}\\
&=
                      \left(
                        \dfrac{1}
                        {
                       1-\dfrac{\min\left\{k_1, 4l_1\right\}}{2(K+L)}
                        }
                        \right)^2-1.
                     \end{split}
           \end{equation}
{\ }

            Observe that,
 by induction, the relation \eqref{normctrl}
            can lead to the same result as \eqref{nsup} for $M$.
            Thus, by comparing \eqref{g2s1} and \eqref{g2s2}
            with \eqref{thm3tanphi}, \eqref{thm4tanphi} and \eqref{Qnice},
            we can include $g=2$ 
            in
            Theorem 4
            under the same dimension assumption
            with no restriction on $l_1$.
            
            
            As a conclusion remark of this section, we point out that,
            if $M=M_1^{k_1}\times N_1^{l_1}\times\cdots\times N_s^{l_s}$
            of $\dim(M)=8$,
            the condition $k_1\geq 4$ is sufficient
            for $C(M)$ to be minimizing in the corresponding Euclidean space. 

{\ }
            
\section{Cones over products of minimal isoparametric hypersurfaces and mixed type}\label{Section5}

In this section we shall consider the case of products of minimal isoparametric hypersurfaces
following the ideas in \S\ref{Section4}
and Theorem 5 will be proved.

An important result of M\"unzner \cite{Munzner}
states that,
in our notations in \S\ref{Section2} for an isoparametric hypersurface of dimension $d$,
 \begin{equation}\label{thetas}
                 \theta_\alpha=\theta_1+\frac{\alpha-1}{g}\pi,  \text{\ \ \ for } 1\leq \alpha\leq g.    
\end{equation}
{Consequently, it follows that (for example, see \cite{TY})
for the unit normal vector field $\xi$} 
{towards $M_+$}
      \begin{equation}\label{MCH}
   d\cdot H=\sum_{\alpha=1}^g m_\alpha \cot \theta_\alpha 
   =\begin{cases}
   m_1g\cot(g\theta_1) &\text{for } g\text{ odd,}\\
   \frac{m_1g}{2}\cot\frac{g\theta_1}{2}
   -\frac{m_2g}{2}\tan\frac{g\theta_1}{2}
   &\text{for } g\text{ even,}
   \end{cases}             
\end{equation}
where $H$ represents the mean curvature.
So for a minimal isoparametric hypersurface,
      \begin{equation}\label{theta1}
  {\ \ \  \ } \theta_1
   =\begin{cases}
   \frac{\pi}{2} &{\ \ \ \ \ \ }g=1,\\
 \arctan\sqrt{\frac{m_1}{m_2}} &{\ \ \ \ \ \ }g=2,\\
 \frac{\pi}{6} &{\ \ \ \ \ \ }g=3,\\
 \frac{1}{2}\arctan\sqrt{\frac{m_1}{m_2}} &{\ \ \ \ \ \ }g=4,\\
\frac{\pi}{12} &{\ \ \ \ \ \ }g=6.
   \end{cases}             
\end{equation}
Hence, the normal radius 
          \begin{equation}\label{NRHyp}
  {\ \ \ \ \ } \phi
   =\begin{cases}
   {\pi} &g=1,\\
 \min\{2\theta_1, \pi-2\theta_1\} &g=2,\\
 \frac{\pi}{3} &g=3,\\
 \min\{2\theta_1, \frac{\pi}{2}-2\theta_1\} &g=4,\\
\frac{\pi}{6} &g=6.
   \end{cases}             
\end{equation}
    Alternatively, 
         \begin{equation}\label{cosphi}
         \,
   \cos\phi
   =\begin{cases}
  0 &\ g=1,\\
1-\frac{2\min\{m_1,\;m_2\}}{m_1+m_2}  &\ g=2,\\
 \frac{1}{2} &\ g=3,\\
 \sqrt{1-\frac{\min\{m_1,\;m_2\}}{m_1+m_2}} &\ g=4,\\
\frac{\sqrt 3}{2} &\ g=6.
   \end{cases}             
\end{equation}
     By $\sqrt{1-\frac{\min\{m_1,\;m_2\}}{m_1+m_2}}<1-\frac{1}{2}\frac{\min\{m_1, \;m_2\}}{m_1+m_2}$,
     it can be concluded that
              \begin{equation}\label{UBcosphi}
   \cos\phi
  \leq
   1-\frac{1}{2d}.       
\end{equation}
     
     Next, we compute the pointwise norm square $\alpha^2$ of the second fundamental form.
     \begin{eqnarray*}
g=1, &\alpha^2=&0,\nonumber\\
g=2, &\alpha^2=&m_1\cot^2\theta_1+m_2\tan^2\theta_1=m_2+m_1=d,\nonumber\\
g=3, &\alpha^2=&m_1(\cot^2\frac{\pi}{6}+\cot^2\frac{3\pi}{6}+\cot^2\frac{5\pi}{6})=6m_1=2d,\nonumber\\
g=4, &\alpha^2=&
m_1(\cot^2\theta_1+\tan^2\theta_1)+m_2(\cot^2(\theta_1+\frac{\pi}{4})+\tan^2(\theta_1+\frac{\pi}{4}))\nonumber\\
&{\ \ \ }=&
m_1\left(\frac{1-2\cos^2\theta_1\sin^2\theta_1}{\cos^2\theta_1\sin^2\theta_1}\right)\\
&{\ \ \ }\ & {\ \ \ \ \ \ \ \ \ \ \ \ \ \ \ \ \ \ \ \ }+
m_2\left(\frac{1-2\cos^2(\theta_1+\frac{\pi}{4})\sin^2(\theta_1+\frac{\pi}{4})}{\cos^2(\theta_1+\frac{\pi}{4})\sin^2(\theta_1+\frac{\pi}{4})}\right)\nonumber\\
&{\ \ \ }=&
m_1\left(\frac{4}{\sin^2 2\theta_1}-2\right)+
m_2\left(\frac{4}{\cos^2 2\theta_1}-2\right)\nonumber\\
&{\ \ \ }=&
m_1\left(\frac{2m_1+4m_2}{m_1}\right)+
m_2\left(\frac{4m_1+2m_2}{m_2}\right)\\
&{\ \ \ }=&6(m_1+m_2)=3d,\nonumber\\
g=6, &\alpha^2=&m_1\left(\sum_{i=0}^{5}\cot^2(\frac{\pi}{12}+\frac{i\pi}{6})\right)=30m_1=5d.\nonumber
\end{eqnarray*}
In summary, 
                 \begin{equation}\label{UBalpha2}
                 \alpha^2=(g-1)d\leq 5d.
                 \end{equation}
                 
                 Assume $M_1^{k_1},\cdots, M_s^{k_s}$ are minimal isoparametric hypersurfaces.
                 Let $\hat \phi$ be the normal radius for the minimal product of $M_1\times\cdots\times M_s$,
                 $\hat \alpha$ the pointwise maximal norm of the second fundamental form in unit normals,
                 and $\hat S=\sum_{i=1}^sk_i$.
                 Then, 
                 by \eqref{normctrl} and \eqref{UBalpha2}, 
                 \begin{equation}\label{UBhatalpha2}
                 \hat \alpha^2\leq 5\hat S.
                 \end{equation}
                 We follow the idea in proving Theorem 3.
                 If the vanishing angle $\theta_F(\hat S+1, \sqrt{5\hat S})$ exists,
                 then similar to \eqref{4c} we have
                  \begin{equation}
                   \tan \theta_0
                  \leq
                  \tan \theta_F(\hat S+1,\sqrt{5\hat S})
                  <
                  \frac{1}{\sqrt{5\hat S}}\sqrt{\dfrac{\hat S}{\hat S-1}}=\dfrac{1}{\sqrt{5(\hat S-1)}},
                 \end{equation}
                 and hence
                  \begin{equation}
                   \tan 2\theta_0
                   < \dfrac{2\cdot\dfrac{1}{\sqrt{5(\hat S-1)}}}{1-\dfrac{1}{{5(\hat S-1)}}}
                 = \dfrac{2\sqrt{5(\hat S-1)}}{5\hat S-6}.
                 \end{equation}
            By the arguments in \S\ref{Section4} and \eqref{UBcosphi} (compared with \eqref{nice}), we get
                  \begin{equation}
                   \tan^2\hat\phi
                   \geq
 \left(\dfrac{1}{1-\dfrac{1}{2\hat S}}\right)^2-1=\dfrac{4\hat S-1}{(2\hat S-1)^2}.
                 \end{equation}
                 It is not hard to check that for $\hat S\geq 5$, $\tan^2\hat\phi>\tan^22\theta_0$.
              
              {\ }   
              
                 Now, the question becomes when $\theta_F$ exists?
                 According to \eqref{BiLi} and \eqref{mono}, we focus on
                 \begin{equation}\label{512}
                \left(\dfrac{\sqrt{19} (\hat S+1)}{12}\right)^2\geq 5\hat S,
                 \end{equation}
                 with positive solutions $\hat S\geq 36$
                 which ensure the existence of $\theta_c$ (and hence that of $\theta_F$).
                 Therefore,
                 cones, of dimension no less than 37, over minimal products of minimal isoparametric hypersurfaces
                 are minimizing.
                 
                 By our arguments and the same relations exhibited in \eqref{UBcosphi} and \eqref{UBhatalpha2},
                 it is clear that 
                 cones, of dimension no less than $37$, over minimal products among minimal isoparametric hypersurfaces
                 and focal submanifolds are minimizing.

                {\ }
                 
                 {\it Remark.} The number $37$ can be improved by more careful calculations.
                 For instance,
                 see the proof of Theorem 4.
                 Also, one can reduce the number by restriction to a subset of minimal isoparametric hypersurfaces and focal submanifolds.
                 There would be
                 several combinations and
                 certain refinements according to the upper bounds of maximal norm squares of second fundamental forms in unit normals.
               
                 {\ }
                 
\appendix
\section*{Appendix}
In \S\ref{NR} we figure out four quantity candidates to be a lower bound of $\tan^2\phi$.
To make arguments in \S\ref{CP}, \S\ref{added} and \S\ref{g2} simple,
a proof of the following useful inequality is 
given here.
\begin{prop*}
For $a,\ b>0$ and $a\neq b$, we have
\begin{eqnarray*}
      \dfrac{4ab}{(a-b)^2}
      &>&
      \min
\left\{
 \left(
                  \dfrac{a+b}{b}
                  \right)^2-1,\
\left(
                  \dfrac{a+b}{a}
                  \right)^2-1
\right\}\\
&=& \min
 \left\{
 \left(
                  \dfrac{1}{1-\frac{a}{a+b}}
                  \right)^2-1,\
\left(
                  \dfrac{1}{1-\frac{b}{a+b}}
                  \right)^2-1
\right\}.
\end{eqnarray*}
\end{prop*}
{\ }\\
{\it Proof. } Assume $b>a>0$.
                One has 
                $$4b^3>3b^3>(b-a)^2(2b+a).$$
                Hence,
                 $$\dfrac{4b}{(a-b)^2}a>\dfrac{2b+a}{b^2}a=\left(\dfrac{a+b}{b}
                  \right)^2-1.$$
                  Similarly, for $a>b>0$, it follows that
                                    $$\dfrac{4a}{(a-b)^2}b>\dfrac{2a+b}{a^2}b=\left(\dfrac{a+b}{a}
                  \right)^2-1.$$
 \hfill$\Box$
\hspace{0.2cm}
\begin{rem}\label{stronger}
For integers $p,q\geq 1$ and $p\neq q$, the above inequality can be improved to
\begin{eqnarray*}
    \dfrac{4pq}{(p-q)^2}
      &\geq&
      \min
\left\{
 \left(
                  \dfrac{p+q}{q-1}
                  \right)^2-1,\
\left(
                  \dfrac{p+q}{p-1}
                  \right)^2-1
\right\}\\
&=& \min
 \left\{
 \left(
                  \dfrac{1}{1-\frac{p+1}{p+q}}
                  \right)^2-1,\
\left(
                  \dfrac{1}{1-\frac{q+1}{p+q}}
                  \right)^2-1
\right\}.
\end{eqnarray*}
Here we think of a nonzero number divided by zero as infinity.
More precisely, when $p>q\geq 1$,
      $$  \dfrac{4pq}{(p-q)^2}
      \geq 
\left(
                  \dfrac{p+q}{p-1}
                  \right)^2-1
=
\left(
                  \dfrac{1}{1-\frac{q+1}{p+q}}
                  \right)^2-1,
$$
with equality if and only if $q=1$.
\end{rem}
\vspace{0.5cm}
\textbf{Acknowledgments}.
We thank Professor Gary Lawlor {for permission of using his graph and table from \cite{Law},}
Professors Michael Kerckhove, Chiakuei Peng, Yuanlong Xin and Doctor Chao Qian for
their interests and useful comments,
{and the referees for kind suggestions.}
The second author wishes to sincerely thank Professor
Leon Simon for many helpful explanations about \cite{HS}, 
Professor H. Blaine Lawson, Jr. for
drawing our attention to \cite{Law},
{and Max-Planck Institute for Mathematics at Bonn for hospitality.}
This work was sponsored in part by {the}
NSFC (Grant Nos. {11331002, 11871282}, {11526048, 11601071}),
Nankai Zhide Foundation,
the Fundamental Research Funds for the Central Universities
(Grant No. 2412016KJ002), 
{the SRF for ROCS, SEM}
 and a Start-up Research Fund from Tongji University. 
{\ }

{\ }
\end{document}